\documentclass{amsart}

\usepackage{url}
\usepackage{hyperref}
\usepackage{verbatim}
\usepackage{amssymb}
\usepackage{color, mathdots}

\newtheorem{theorem}{Theorem}[section]
\newtheorem{proposition}[theorem]{Proposition}
\newtheorem{lemma}[theorem]{Lemma}
\newtheorem{corollary}[theorem]{Corollary}

\theoremstyle{definition}
\newtheorem{definition}[theorem]{Definition}
\newtheorem{conjecture}[theorem]{Conjecture}

\theoremstyle{remark}
\newtheorem{remark}[theorem]{Remark}
\newtheorem{example}[theorem]{Example}
\newtheorem{question}[theorem]{Question}

\def \gg {\mathfrak g}
\def \LL {\mathcal L}
\def \MM {\mathfrak M}

\def \SS {{\bf S}}
\def \HH {{\bf H}}
\def \KK {{\bf K}}
\def \C {\mathfrak C}
\def \CP {\mathfrak C_{\text{Poi}}}

\numberwithin{equation}{section}
 
\newcommand{\mf}{\mathfrak}
\newcommand{\mb}{\mathbb}
\newcommand{\mc}{\mathcal}
 
 \newcommand{\NN}{\mb N}
 \newcommand{\ZZ}{\mb Z}
\newcommand{\DMO}{\DeclareMathOperator}
\DMO{\height}{ht}
\DMO{\ad}{ad}
\newcommand{\BB}{\mf B}
\newcommand{\wh}{\widehat}
\newcommand{\ff}{\footnote}

\newcommand{\beq}{\begin{equation}}
\newcommand{\eeq}{\end{equation}}
\newcommand{\del}{\partial}

\newcommand{\kk}{{\Bbbk}}

\def \WW {{\bf W}}
\newcommand{\Vir}{V\!ir}
\DMO{\alg}{alg}

\title[A Poisson basis theorem for symmetric algebras]{A Poisson basis theorem for symmetric algebras of infinite-dimensional Lie algebras}

\author{Omar Le\'on S\'anchez}
\address{Omar Le\'on S\'anchez\\
University of Manchester\\
Department of Mathematics\\
Manchester, M13 9PL, U.K.}
\email{omar.sanchez@manchester.ac.uk}
\thanks{{\em Acknowledgements}: The first author was partially supported by EPSRC grant EP/V03619X/1}

\author{Susan J. Sierra}
\address{Susan J. Sierra\\
University of Edinburgh\\
School of Mathematics\\
Edinburgh EH9 3JZ, U.K.}
\email{s.sierra@ed.ac.uk}

\date{\today}

\begin{document}

\begin{abstract}
We consider when the symmetric algebra of an infinite-dimensional Lie algebra, equipped with the natural Poisson bracket, satisfies the ascending chain condition (ACC) on Poisson ideals. 
We define a combinatorial condition on a graded Lie algebra which we call \emph{Dicksonian} because it is related to Dickson's lemma on finite subsets of $\NN^k$.  
Our main result is:

\medskip
\noindent {\bf Theorem.  } If $\gg$ is  a Dicksonian graded Lie algebra over a field of characteristic zero, then the symmetric algebra $S(\gg)$ satisfies the ACC on radical Poisson ideals. 
\medskip

As an application, we establish this ACC for the symmetric algebra of  any graded simple Lie algebra of polynomial growth over an algebraically closed field of characteristic zero, and for the symmetric algebra of the Virasoro algebra. We also derive some consequences connected to the Poisson primitive spectrum of finitely Poisson-generated algebras.
\end{abstract}

\maketitle

\tableofcontents

\section{Introduction}

{
Let $\kk$ be a field of characteristic zero. Recall that a Poisson algebra $A$ is a commutative $k$-algebra equipped with a skew-symmetric $\kk$-bilinear map 
$$\{-,-\}:A\times A \to A$$ 
such that for each $a\in A$ the map $\{a,-\}:A\to A$ yields a derivation. By a Poisson ideal we mean an ideal $I$ of $A$ such that $\{a,I\}\subseteq I$ for all $a\in A$.} This paper deals with, and is motivated by, Ascending Chain Conditions on Poisson ideals in certain classes of Poisson algebras. More precisely, we study noetherianity of certain systems of Poisson ideals in the symmetric algebra $S(\gg)$ of a $\mathbb Z$-graded Lie $\kk$-algebra $\gg$. 
{
Here $S(\gg)$ can be identified with the polynomial ring $\kk[\MM]$ with $\MM$ a $\kk$-basis for $\gg$ (understood as formal variables) and the canonical Poisson structure is induced from the Lie bracket on $\gg$}. 

{
The problem above} is related  to the following general question: for which Lie algebras $\gg$ is the system of two-sided ideals of the enveloping algebra $U(\gg)$ noetherian? For instance,  in \cite[Conjecture 1.3]{PS} the following was conjectured:

\begin{conjecture}\label{conjecturewitt}
Let $W_+$ be the positive Witt algebra, which has basis $\{e_n : n \in \ZZ_{\geq 1}\}$ with 
\beq \label{Witt}
[e_n, e_m ]= (m-n) e_{n+m}.
\eeq
 The system of two-sided ideals of $U(W_+)$ satisfies the ACC. 
\end{conjecture}

Via the associated graded construction and using \cite[Proposition 1.6.8]{MR}, Conjecture~\ref{conjecturewitt} would follow if we knew that the symmetric algebra $S(W_+)$ equipped with its natural Poisson structure had the ACC on Poisson ideals. While this remains open, in \cite[Corollary 2.17]{PS} it is established that $U(W_+)$ satisfies the ACC on (two-sided) ideals whose associated graded Poisson ideal in $S(W_+)$ is radical. This is a consequence of \cite[Theorem 2.15]{PS} which states that the symmetric algebra $S(W_+)$  satisfies the ACC on radical Poisson ideals. 

As pointed out in \cite[Remark 2.18]{PS}, techniques from differential algebra can be useful in the study of the Poisson ideal structure of $S(W_+)$. In this paper we exploit this idea to prove a general Poisson basis theorem for symmetric algebras $S(\gg)$ where the Lie algebra $\gg$ is graded 
{
(namely, $\mathbb Z$-graded with finite-dimensional homogeneous components)} and satisfies a certain  ``combinatorial noetherianity" condition, which we call {\em Dicksonian} as it relates to Dickson's lemma on subsets of $\NN^k$. More precisely, taking our cue from techniques involved in the differential basis theorem of Kolchin \cite[Chapter III, \S4]{kolchin},  in Section \ref{eliminationalg} we prove our main result:

\begin{theorem}\label{basisPoisson}
{\rm (Theorem~\ref{thm:basisPoissonrestate})}
Let $\kk$ be a field of characteristic zero and $\gg$ a graded Lie $\kk$-algebra. If $\gg$ 
is Dicksonian,
then the Poisson algebra $S(\gg)$ has ACC on radical Poisson ideals.
\end{theorem}

In Section \ref{leading} we explain what we mean by a Lie algebra 
being {\em Dicksonian}. 
We also provide, in Lemma~\ref{exconditions}, sufficient conditions that guarantee
this property for $\mf g$.
 One easily checks that the positive Witt algebra $W_+$ satisfies these conditions, and thus $S(W_+)$ has ACC on radical Poisson ideals, as already pointed out in \cite[Theorem 2.15]{PS}. One can also easily check that the full Witt algebra $W$ (which has basis $\{ e_n : n \in \ZZ\}$ and Lie bracket \eqref{Witt}) satisfies these \ conditions, and thus the symmetric algebra $S(W)$ also has ACC on radical Poisson ideals.
 To our knowledge this does not appear elsewhere.

\begin{remark}\label{rem1.3}
We note that the conclusion of Theorem \ref{basisPoisson} cannot generally be strengthened to the ACC on the whole system of Poisson ideals. Consider the following example. Let $D$ be the Lie algebra generated as a $\kk$-vector space by $x_1,x_2,\dots,$ and $y$, with bracket 
$$[x_i,x_j]=0 \quad \text{and} \quad [y,x_i]=x_{i+1}.$$
We can equip $D$ with the $\mathbb Z$-grading where $x_i$ has degree $i$ and $y$ has degree $1$. One readily checks that $D$ has ACC on Lie ideals and, moreover, it is Dicksonian. Thus $S(D)$ has ACC on radical Poisson ideals. However, $S(D)$ does not have ACC on arbitrary Poisson ideals. For example, the chain 
$$[x_1^2]\; \subseteq \; [x_1^2, x_2^2]\; \subseteq \; [x_1^2,x_2^2,x_3^2]\; \subseteq \;  \cdots$$
is strictly increasing. Here $[S]$ denotes the Poisson ideal generated by $S\subseteq D$.

We note that in $S(W_+)$ the chain $[e_1^2] \subseteq [e_1^2, e_2^2] \subseteq [e_1^2,e_2^2,e_3^2]\subseteq\cdots $ does stabilise, by \cite[Corollary~4.8]{PS}.  This illustrates the delicacy of the noetherianity questions we discuss here.
\end{remark}

\smallskip

We further note that the conclusion of Theorem~\ref{basisPoisson} cannot be achieved without assuming some suitable chain condition on the ideals of the Lie algebra $\gg$.

\begin{lemma}
Let $\kk$ be a field and let $\gg$ be a Lie $\kk$-algebra. If $\gg$ is not noetherian, then $S(\gg)$ does {\bf not} have ACC on radical Poisson ideals.
\end{lemma}
\begin{proof}
As $\gg$ is not noetherian, there exists a strictly increasing chain of Lie-ideals $I_1\subset I_2\subset \cdots$. Let $P_i$ be the ideal of $S(\gg)$ generated by $I_i$ for $i=1,2,\dots$. Then, as $I_i$ is a Lie-ideal, it can easily be checked that $P_i$ is a Poisson ideal of $S(\gg)$. Also, as the $P_i$'s are generated by linear terms, they are all prime. Finally, one readily checks that the chain $P_1\subset P_2\subset\cdots $ is strictly increasing (indeed any element in $I_{i+1}\setminus I_i$ is not in $P_i$).
\end{proof}

Examples of non-noetherian Lie $\kk$-algebras include free Lie algebras over $\kk$ in at least two generators. For instance, in \cite[Theorem 1.1]{AR96}, it is shown that if $I$ is a nontrivial Lie-ideal of a free Lie algebra then $[I,I]$ is not finitely generated as a Lie-ideal. These examples show that even when the Lie algebra $\gg$ is finitely generated it is not generally the case that $S(\gg)$ has ACC on radical Poisson ideals.

\medskip

Our long-term goal is to apply Theorem~\ref{basisPoisson} to a wide class of Poisson algebras: for instance, to all symmetric algebras of noetherian graded Lie algebras. 
We are thus far not aware of a counter-example, although we caution that noetherian Lie algebras can be quite wild; see \cite{Greenfeld}.  
A somewhat more accessible short-term goal is to restrict ourselves to Lie algebras with well-behaved growth. In fact, we expect:

\begin{conjecture}\label{vision}
Let $\kk$ be an algebraically closed  field of characteristic zero and $\gg$ a graded Lie $\kk$-algebra of polynomial growth (also called finite growth). If $\gg$ has ACC on Lie ideals, then $\gg$ 
is Dicksonian. 
(As a consequence of Theorem~\ref{basisPoisson}, $S(\gg)$ would have ACC on radical Poisson ideals). 
\end{conjecture}

A reasonable place to start towards proving Conjecture \ref{vision}, is the case when $\gg$ is a simple graded Lie algebra 
{
(namely, $\gg$ has no non-trivial homogenous ideals)} of polynomial growth. When the field $\kk$ is algebraically closed, these Lie algebras have been classified by Mathieu \cite{mathieu}; they are either finite-dimensional, loop algebras, Cartan algebras, or the Witt algebra. In Section \ref{examples}, we verify the conjecture for all these classes of Lie algebras, and then prove:

\begin{corollary}\label{simplecase}
{\rm (Corollary~\ref{cor:simplerestate})}
Let $\kk$ be a field of characteristic zero  and $\gg$ a Lie $\kk$-algebra. If $\kk^{\text{alg}}\otimes_\kk\gg$ is a simple graded Lie $\kk^{\text{alg}}$-algebra of polynomial growth, then the symmetric algebra $S(\gg)$ has ACC on radical Poisson ideals. 
\end{corollary}

We also prove that the  Virasoro algebra and several related Lie algebras are Dicksonian, and thus their symmetric algebras have  ACC on radical Poisson ideals.

In our final section, Section \ref{PoiDME}, we consider the Poisson version of the Dixmier-Moeglin equivalence relating primitive, rational, and locally closed prime ideals of enveloping algebras of finite-dimensional Lie algebras. We make several remarks on the Poisson primitive spectrum of Poisson algebras of countable vector space dimension; in particular, those that are finitely Poisson-generated. For instance, we prove in Theorem \ref{implications}, 
{
under the assumption that $\kk$ is uncountable}, that the notions of Poisson-primitive and Poisson-rational coincide in this setting. We then pay a closer look at $S(W_+)$ in this context. 
{
We remind the reader that a prime Poisson ideal $P$ of a Poisson algebra $A$ is said to be Poisson-locally closed if $A$ is a locally closed point in the Poisson spectrum of $A$ (with respect to the induced Zariski-topology); $P$ is called Poisson-primitive if it is the Poisson core of a maximal ideal of $A$; and, finally, $P$ is called Poisson-rational if the Poisson center of the fraction field of $A/P$ is algebraic over $\kk$.}

\medskip

{\bf Acknowledgements: } We thank Rekha Biswal and  Alexey Petukhov for useful discussions on twisted loop algebras and Poisson primitive ideals, respectively. 
We also thank the anonymous referee for carefully reading a previous version and for their several helpful comments and suggestions.

\medskip

{\bf Notation:} For us, $\NN = \{0, 1, 2, \dots \}$.

\

\section{A review on radical conservative systems}

In this section we review some standard results on radical divisible conservative systems in arbitrary commutative rings. We follow closely \cite[Chapter 0, \S7 - \S9]{kolchin} (where the omitted proofs appear). We then specialize these to the context of Poisson algebras where the main result for us is Theorem~\ref{conservativePoisson}.

We fix a commutative ring $R$ with unit. Recall that a conservative system $\C$ of $R$ is a collection of ideals of $R$ with the following two properties:
\begin{enumerate}
\item [(CS1)] the intersection of elements in $\C$ is in $\C$, and
\item [(CS2)] the union of a chain (totally ordered set by inclusion) from $\C$ is again in $\C$.
\end{enumerate}

For example, the collection of all ideals of $R$ is a conservative system; as is, on the other extreme, the collection consisting just of $R$. 

Let $\C$ be a conservative system of $R$. We refer to the elements of $\C$ as $\C$-ideals.  Given an arbitrary subset $A$ of $R$, we denote by $(A)_\C$ the intersection of all the $\C$-ideals containing $A$. Note that $(A)_\C$ is in $\C$ by condition (CS1). Hence, we call it the {\em $\C$-ideal $\C$-generated by $A$}. If a $\C$-ideal $I$ is of the form $(\Sigma)_\C$ for some finite set $\Sigma$, we say that $I$ is {\em finitely $\C$-generated}.

Recall that given an ideal $I$ of $R$ and $s\in R$, the {\em division of $I$ by $s$} is the ideal $I:s$ defined by $\{r\in R: rs\in I\}$. 


\begin{definition} \
\begin{enumerate} 
\item [(i)] The conservative system $\C$ is called {\em divisible} if for all $I\in \C$ and $s\in R$, we have $I:s\in \C$. 
\item [(ii)] The conservative system $\C$ is called {\em radical } if all of its elements are radical ideals. 
\item [(iii)] A conservative system $\C$ is said to be {\em noetherian} if it satisfies the ACC on $\C$-ideals (equivalently, every $\C$-ideal is finitely $\C$-generated).
 \end{enumerate}
\end{definition} 

\begin{example}
The collection of all radical ideals of $R$ is a radical divisible conservative system.
\end{example}

The following two lemmas are two of the main ingredients for the proof of Theorem~\ref{generalbasis} below.

\begin{lemma}\label{nice1}\cite[\S0.8 Lemma 7]{kolchin}
Suppose $\C$ is a radical divisible conservative system of $R$. If $T$ and $S$ are arbitrary subsets of $R$, then
$$(T\cdot S)_\C=(T)_\C\cap (S)_\C.$$
Here $T\cdot S$ denotes the set of products of the form $ts\in R$ with $t\in T$ and $s\in S$.
\end{lemma}

\begin{lemma}\label{nice2}\cite[\S0.9 Lemma 8]{kolchin}
Suppose $\C$ is a radical divisible conservative system of $R$. If $\C$ is not noetherian, then there is a $\C$-ideal that is maximal (with respect to inclusion) among the $\C$-ideals that are not finitely $\C$-generated and, more importantly, any such $\C$-ideal is prime. 
\end{lemma}

We now prove an algebraic ``basis theorem", which we will use in the Poisson context later on.

\begin{theorem}\label{generalbasis}
Let $\C$ be a radical divisible conservative system of $R$. Assume that
\begin{enumerate}
\item [(*)] if $P$ is a prime $\C$-ideal, then there is a finite $\Sigma\subset P$ and $s\in R\setminus P$ such that $P=(\Sigma)_\C:s$.
\end{enumerate}
Then, $\C$ is noetherian.
\end{theorem}
\begin{proof}
The proof can be deduced from arguments in \cite[Chapter 0, \S9]{kolchin}, but as it is not explicitly stated there, we prove it here. 

Towards a contradiction, assume $\C$ is not noetherian. Then, by Lemma~\ref{nice2}, there is a maximal $\C$-ideal $M$, with respect to inclusion, among the $\C$-ideals that are not finitely $\C$-generated, and $M$ is prime. By assumption (*), there is a finite set $\Sigma\subset M$ and $s\in R\setminus M$ such that $M=(\Sigma)_\C: s$. It follows that $s\cdot M\subseteq (\Sigma)_\C$, Furthermore, as $s\notin M$, by choice of $M$ there is a finite $\Phi\subset M$ such that $(s, M)_\C=(s,\Phi)_\C$. Thus, using Lemma~\ref{nice1}, we get
$$M=M\cap (s,M)_\C=M\cap(s,\Phi)_\C=(s\cdot M, \Phi)_\C=(\Sigma, \Phi)_\C.$$
This contradicts the fact that $M$ is not finitely $\C$-generated.
\end{proof}

We conclude this review on conservative systems with a decomposition-type theorem for radical ideals. If $\C$ is a conservative system of $R$ and $I$ is a radical $\C$-ideal, by a $\C$-component of $I$ we mean a minimal element of the set, ordered by inclusion, of prime $\C$-ideals that contain $I$. We have the following

\begin{proposition}\cite[\S0.8 Proposition 1, \S0.9 Theorem 1]{kolchin}
Assume $\C$ is a radical divisible conservative system of $R$ and let $I$ be a $\C$-ideal. Then, the following hold:
\begin{enumerate}
\item [(i)] $I$ is the intersection of its $\C$-components,
\item [(ii)] If $I$ is the intersection of finitely many prime ideals none of which contains the other, then these prime ideals are in $\C$ and are the $\C$-components of $I$.
\item [(iii)] If $\C$ is noetherian, $I$ is the intersection of a finite set of prime $\C$-ideals none of which contains the other. This finite set is unique, being precisely the set of $\C$-components of $I$.
\end{enumerate}
\end{proposition}

\subsection{Applications to radical Poisson ideals}\label{specialPoi}

We now specialize some of the results above to Poisson algebras (the rest of the translations are left to the interested reader). Let $(A,\{-,-\})$ be a Poisson algebra over a field $\kk$; note we make no assumption on the characteristic at this point. Let $\CP$ be the collection of all radical Poisson ideals of $A$. Then $\CP$ is a radical divisible conservative system of $A$. Indeed, all conditions are more or less clearly satisfied; we only provide details on divisibility. Let $I\in \CP$ and $s\in A$. We prove that $I:s\in \CP$ (i.e., $I:s$ is a radical Poisson ideal of $A$). Radicality easily follows. Now let $g\in A$, and suppose $f\in I:s$. Then $sf$ and $s^2f$ are in $I$; and so $\{g,s^2f\}\in I$ (as $I$ is Poisson). But $\{g,s^2f\}=\{g,s\}2sf+s^2\{g,f\}$, and so $s^2\{g,f\}\in I$. Hence, $\{g,f\}\in I:s$ (as $I$ is radical), and so $I:s$ is Poisson.

For $\Sigma$ any subset of $A$, we let $\{\Sigma\}$ denote the radical Poisson ideal of $A$ generated by $\Sigma$. Note that $\{\Sigma\}=(\Sigma)_{\CP}$. Here are the relevant specializations:

\begin{theorem}\label{conservativePoisson}
Let $(A,\{-,-\})$ be a Poisson algebra over a field $\kk$, and let $\CP$ be the system of radical Poisson ideals.\begin{enumerate}
\item [(i)] Assume that for every prime Poisson ideal $P$ of $A$ there is a finite set $\Sigma\subset P$ and $s\in A\setminus P$ such that $P=\{\Sigma\}:s$. Then, $\CP$ is noetherian. 
\item [(ii)] Let $I$ be a $\CP$-ideal (i.e., a radical Poisson ideal). Then, $I$ is the intersection of prime Poisson ideals. Furthermore, if the system $\CP$ is noetherian, then $I$ is the intersection of a finite set of prime Poisson ideals none of which contains the other (this set is unique and its elements are called the Poisson-components of $I$).  \qed
\end{enumerate}
\end{theorem}

Below we will use part (i) of this theorem to prove noetherianity of the system $\CP$ for Poisson algebras of the form $S(\gg)$ where $\gg$ is a Dicksonian graded Lie algebra over a field $\kk$ of characteristic zero.

\section{Dicksonian  Lie algebras}\label{leading}

Our goal in this section is to define the combinatorial condition used in Theorem~\ref{basisPoisson}, which is named for Dickson's lemma on subsets of $\NN^k$ (see Theorem~\ref{Dickson}).  We consider orderings on a Lie algebra $\mf g$ and the information they give us on Lie ideals of $\mf g$.  We will then define the concept of a {\em leading-Dicksonian sequence}:  a sequence of (pairs of) elements of $\mf g$ satisfying a certain chain condition.   A Lie algebra is {\em Dicksonian} if it has no infinite leading-Dicksonian sequence.   We will see in  Section~\ref{eliminationalg} that Dicksonian Lie algebras have an elimination algorithm, which allows us to derive striking consequences for radical Poisson ideals of their symmetric algebras.

Throughout this section we assume that $\gg$ is a ($\mathbb Z$-)graded Lie algebra over a field $\kk$ of characteristic zero. Namely, $\gg$ is a Lie $\kk$-algebra equipped with a decomposition $\gg=\bigoplus_{n\in \mathbb Z}\gg_n$ such that $[\gg_n,\gg_m]\subseteq \gg_{n+m}$ and the homogeneous components $\gg_n$ have finite dimension. We refer the reader to \cite[\S1]{mathieu} for basic facts on graded Lie algebras.

Let $\MM$ be a $\kk$-basis of $\gg$ consisting of homogeneous elements equipped with a total order $(\mathfrak M,<)$ that is compatible with the grading (i.e., larger in degree implies larger in the order $<$). 
As the $\mf g_n$ are finite-dimensional, $\MM$ has the order type of a subset of $\ZZ$.  
We let $\MM_+$ and $\MM_-$ denote the elements of $\MM$ of positive and negative degree (with respect to the grading of $\gg$), respectively. Furthermore, for nonzero $e\in \gg$, we let $\ell_{+}(e)$, respectively $\ell_{-}(e)$, denote the largest, respectively smallest, element of $\MM$ with respect to $<$ that appears in $e$ when written as a $\kk$-linear combination of elements of $\MM$.
 We call $\ell_{+}(e)$ the upper-leader and $\ell_{-}(e)$ the lower-leader of $e$. As convention, we set $\ell_{\pm}(0)=0$. 

\begin{definition}\label{Di}
Let $\mathcal I_{\pm}$ be the set of (nonempty) finite tuples of elements from $\MM_{\pm}$. For any ${\bf i}=(M_1,\dots,M_n)\in \mathcal I_{\pm}$ we let $D_{\bf i}^{\pm}$ be the operator on $\MM$ given by
$$M\mapsto D_{\bf i}^{\pm}(M):=\ell_{\pm}([[[M,M_1],M_2],\cdots,M_n]).$$
\end{definition}

For $M\in \MM$, we set $\LL_+(M)$ to be the set of elements of $\MM$ of the form
$$D_{\bf i}^+(M),$$
for ${\bf i}\in \mathcal I_+$, such that $D_{\bf i}^+(N)<D_{\bf i}^+(M)$ for all $N<M$ whenever $D_{\bf i}^+(N)$ is nonzero. Similarly, $\LL_-(M)$ denotes the set of elements of $\MM$ of the form
$$D_{\bf i}^-(M),$$
for ${\bf i}\in\mathcal I_-$, such that $D_{\bf i}^-(N)>D_{\bf i}^-(M)$ for all $N>M$ whenever $D_{\bf i}^-(N)$ is nonzero. 

We summarise the notation in Table~\ref{table1}.
\begin{table}[h]
\begin{center}
    \caption{Notation for operations on an ordered basis $\mf M$ of $\mf g$}
    \label{table1}
\begin{tabular}{cc}     
    $\mf M_{\pm}$ & elements of  $\mf M$ of positive (negative) degree \\
    $\ell_{\pm}(e)$ & largest (smallest) element of $\mf M$ occurring in $e \in \mf g$ \\
    $\mc I_{\pm}$ & finite tuples from $\mf M_{\pm}$ \\
    $D_{\bf i}^{\pm}(M)$ & $\ell_{\pm}([[[M,M_1],M_2],\cdots,M_n])$ where ${\bf i}=(M_1,\dots,M_n)\in \mathcal I_{\pm}$ \\
    $\mc L_+(M) $ & $\{ D_{\bf i}^{+}(M) : D_{\bf i}^{+}(N) < D_{\bf i}^{+}(M) \mbox{ if ${\bf i} \in \mathcal I_+, N < M, D_{\bf i}^+(N) \neq 0$ } \} $\\
    $\mc L_-(M) $ & similar
        \end{tabular}
  \end{center}
  \end{table}
  
\begin{definition}\
\begin{enumerate}
\item [(i)] A sequence of {\bf distinct} pairs 
$$((M_i,N_i))_{i=1}^{n}$$
from $\MM$ with $n\leq \omega$ and $M_i\leq N_i$ is said to be {\em leading-Dicksonian} (with respect to the order $<$ of $\MM$) if $M_j\notin \LL_-(M_i)$ and $N_j\notin \LL_+(N_i)$ for all $i<j$.
\item [(ii)] We say that $\gg$ 
is {\em Dicksonian} if there is a basis of homogeneous elements with an order compatible with the grading such that, with respect to this order, there is no infinite leading-Dicksonian sequence. 
\end{enumerate}
\end{definition}

\begin{example}\label{Wittranked}
The (full) Witt algebra $W$ is Dicksonian.
Indeed, choosing the standard $\kk$-basis $(e_n:n\in \mathbb Z)$, for which $[e_i,e_j]=(j-i)e_{i+j}$, one sees that the order in $W$ given by 
$$e_i<e_j \quad \text{ iff } \quad  i<j$$ 
has the desired properties. For example, $\LL_+(e_1) = \{e_i : i > 2\}$, and if $n\neq 1$ then $\LL_+(e_n) = \{e_i : i > n\}$.
The positive Witt algebra 
$$W_+=\text{span}_\kk(e_i:i\geq 1),$$ and  the Cartan algebra 
$$\WW_1=\text{span}_\kk(e_i:i\geq -1)$$ are also Dicksonian, by a similar argument.
\end{example}

Extending the above example, we have

\begin{lemma}\label{exconditions}
Suppose $\gg$ is a graded Lie algebra with an ordered basis $(\MM,<)$ where $\MM$ consists of homogeneous elements and the order is compatible with the grading. Assume that the Lie bracket of two basis elements is a scalar multiple of a basis element and that the following condition holds
\begin{enumerate}
\item [($\dagger$)] if $M_1,M_2, M\in \mathfrak M_{\pm}$ are such that $[M_1,M]$ and $[M_2,M]$ are nonzero and $M_1< M_2$, then $\ell_{\pm}([M_1,M])<\ell_{\pm}([M_2,M])$.
\end{enumerate}
If $\gg_+$ and $\gg_-$ have ACC on graded Lie ideals, then $\gg$ has no infinite leading-Dicksonian sequence.
\end{lemma}
Note in $(\dagger)$ that here $[M_1, M]$ and $[M_2, M]$ are scalar multiples of a basis element, and $\ell_{\pm}$ simply extracts this element.
\begin{proof}
Towards a contradiction, let  $((M_i,N_i))_{i=1}^\infty$ be  an infinite leading-Dicksonian sequence. Note that there are either infinitely many $M_i$'s in $\mathfrak M_-$ or infinitely
 {
 many} $N_i$'s in $\mathfrak M_+$. Without loss of generality assume the latter, we will show that $\gg_+$ has a strictly ascending chain of graded Lie ideals, contradicting our assumption.

We thus have an infinite sequence $(N_i)_{i=1}^{\infty}$ of homogenous elements of $\gg_+$ with the property that $N_j\notin \mathcal L_+(N_i)$ for all $i<j$. It suffices to show that this latter condition implies that $N_j$ is not in the Lie ideal generated by $N_1, \dots,N_{j-1}$ in $\gg_+$. Suppose towards a contradiction that $N_j$ is in this Lie ideal. By our assumption that the bracket of basis elements yields a scalar multiple of a  basis element, we must have that $N_j=D^+_{\bf i}(N_i)$ for some $\bf i\in \mathcal I_+$ and $i<j$. This yields, by condition ($\dagger$), that $N_j\in \mathcal L_+(N_i)$, a contradiction. The result follows.
\end{proof}

\begin{example}\label{example1}\
\begin{enumerate}
\item Consider the loop algebra $\wh{\mf{sl}_2} := \mf{sl}_2(\kk)[t, t^{-1}]$.  Letting $e, f, h$ be the standard basis of $\mf{sl}_2(\kk)$, let $\mf M = \{ et^i, ft^j, h t^k : i, j, k \in \ZZ \}$.  Give $e, f, h, t$ degrees $1, -1, 0, 3$, respectively, and order elements of $\mf M$ by degree.  By Lemma~\ref{exconditions}, $\wh{\mf{sl}_2}$ is Dicksonian.
\item It also follows from Lemma~\ref{exconditions} that the  Lie algebra $D$ of Remark~\ref{rem1.3} is Dicksonian.
\end{enumerate}
\end{example}

\begin{remark}\label{rem:newcond}
If $\mf g $ has a basis $\MM$ so that for all $M \in \MM$ we have
\beq\label{cofinite}
\LL_+(M) \cup \LL_-(M)\mbox{ is cofinite in } \MM,
\eeq
then $\mf g$ is easily seen to be Dicksonian; noting that $\MM$ has order type of a subset of $\ZZ$.
This gives an alternate proof that $W$, $W_+$, and $\WW_1$ are Dicksonian.
\end{remark}

\begin{example}
Let $\Vir$ be the {\em Virasoro algebra}, which has  basis $\{e_n : n \in \ZZ\} \cup \{z\}$ and Lie bracket 
\[[e_n, e_m] = (m-n) e_{n+m} + \frac{n^3-n}{12} \delta_{n+m, 0} z, \quad [e_n, z]=0. \]
Ordering the basis by
\[ \dots < e_{-2} < e_{-1} < z < e_0 < e_1 < e_2 < \dots,\]
one sees that \eqref{cofinite} is satisfied for all basis elements except $z$.
Thus by Remark~\ref{rem:newcond}, $\Vir$ is Dicksonian.
\end{example}

\section{Elimination algorithms in $S(\gg)$}\label{eliminationalg}

In this section we prove our main result,  Theorem~\ref{basisPoisson}. One of the key ingredients is an \emph{Elimination Algorithm} result for $S(\gg)$ that we prove in Theorem~\ref{elimination} below. We take our cue/presentation from the elimination theory of differential polynomial rings; see for instance \cite[Chapter I]{kolchin}.

\medskip

\noindent {\bf Assumptions:} Throughout this section we assume that $\gg$ is a graded Lie algebra over a field $\kk$ of characteristic zero  and $\MM$ is a $\kk$-basis consisting of homogeneous elements equipped with an order $<$ compatible with the grading. Define $\MM_+$ and $\MM_-$ as in Section~\ref{leading}.   Recall that the symmetric algebra $S(\gg)$ is the polynomial ring in the formal variables $\MM$ over $\kk$ (i.e., $S(\gg)=\kk[\mathfrak \MM]$) and it carries a natural Poisson bracket, that we denote by $\{-,-\}$, induced from the Lie bracket on $\gg$.

\medskip

Let $f$ be a nonconstant element in $S(\gg)$; i.e., a nonconstant polynomial in $\kk[\mathfrak M]$. We define the upper leader of $f$, denoted by $\ell_{+}(f)$, to be the largest element in $\mathfrak M$ (according to the fixed order $<$) that  appears (nontrivially) in $f$. Then, $f$ can be written in the form
$$f=\sum_{i=0}^{d} g_i (\ell_{+}(f))^ i$$
with the $g_i$'s in $\kk[\mathfrak M]$ having leaders strictly smaller than $\ell_{+}(f)$, and $g_d$ nonzero. We define the upper-degree of $f$, denoted $d_{+,f}$, to be this $d$. We call $g_d$ the upper-initial of $f$, denoted $i_{+,f}$. The upper-separant of $f$ is defined to be 
$$s_{+,f}=\frac{\partial f}{\partial \ell_{+}(f)}.$$
The upper-rank of $f$ is defined as $rk_+(f)=(\ell_{+}(f),d_{+,f})$. We can compare elements from $S(\gg)$ lexicographically by upper-rank. Note that $s_{+,f}$ and $i_{+,f}$  both have lower upper-rank than $f$.

In a similar fashion one defines the lower-leader of $f$, denoted by $\ell_{-}(f)$, as the smallest element of $\MM$ that appears in $f$. The notions of lower-degree $d_{-,f}$, lower initial $i_{-,f}$,  lower-separant $s_{-,f}$, and lower-rank $rk_-(f)$ are defined similarly.

\begin{definition}
Recall that $\mathcal I_{\pm}$ denotes the set of (nonempty) finite tuples of elements from $\MM_{\pm}$. For any ${\bf i}=(M_1,\dots,M_n)\in \mathcal I_{\pm}$ we let $D_{\bf i}$ be the operator on $S(\gg)$ given by
$$f\mapsto D_{\bf i}(f):=\{\{\{f,M_1\},M_2\},\cdots,M_n\}.$$
Note that $D_{\bf i}(f)$ is always in the Poisson ideal $[f]$ generated by $f$.
\end{definition}

\medskip

For $f\in S(\gg)$ we let $s_f=s_{+,f}\cdot s_{-,f}$. For $\Lambda$ a finite subset of $S(\gg)$ we let $\ell_{+}(\Lambda)=\max\{\ell_{+}(f):f\in \Lambda\}$ and $\ell_{-}(\Lambda)=\min\{\ell_{-}(f):f\in \Lambda\}$.

\medskip

We summarise this notation in Table~\ref{table2}.

\begin{table}[h]
\begin{center}
    \caption{Notation for operations on $S(\mf g)$}
    \label{table2}
\begin{tabular}{cc}     
    $\ell_\pm(f)$ & upper- (lower-) leader of $f \in S(\mf g)$ \\
    $d_{\pm, f}$ & uppper- (lower-) degree of $f$ \\
    $i_{\pm, f}$ & upper- (lower-) initial of $f$ \\
    $s_{\pm f}$ & upper- (lower-) separant of $f$, $\frac{\partial f}{\partial \ell_\pm(f)}$ \\
    $rk_{\pm}(f)$ & upper- (lower-) rank of $f$, $(\ell_{\pm}(f), d_{\pm,f})$ \\
    $D_{\bf i}(f) $ & $\{\{\{f,M_1\},M_2\},\cdots,M_n\}$ for ${\bf i} = (M_1, \dots, M_n)\in \mathcal I_{\pm}$\\
    $s_f$ & $ s_{+,f}\cdot s_{-,f}$ \\
    $\ell_{+}(\Lambda)$ & $\max\{\ell_{+}(f):f\in \Lambda\}$ \\
    $\ell_{-}(\Lambda)$ & $\min\{\ell_{-}(f):f\in \Lambda\}$.
        \end{tabular}
  \end{center}
  \end{table}

\medskip

The content of the next lemma is that if the sets $\mc L_{\pm}(e)$ are large enough, one can control the largest/smallest element of $\mf M$ appearing in some $D_{\bf i}(f)$.

\begin{lemma}\label{separant}
Let $f\in S(\gg)$ be nonconstant and ${\bf i}\in \mathcal I_{\pm}$. If $D_{\bf i}^{\pm}(\ell_{\pm}(f))\in \LL_{\pm}(\ell_{\pm}(f))$, then 
$$\ell_{\pm}(D_{\bf i}(f))=D_{\bf i}^{\pm}(\ell_{\pm}(f))$$ 
and 
$$D_{\bf i}(f)=\alpha s_{\pm,f}D_{\bf i}^{\pm}(\ell_{\pm}(f))+h_{\pm}$$
for some $\alpha\in k^*$ and $h_{+}, h_{-}\in S(\gg)$, with 
$$\ell_-(f)\leq \ell_{-}(h_+)\leq \ell_{+}(h_+)<D_{\bf i}^{+}(\ell_{+}(f))\;  \mbox{ for } {\bf i} \in \mc I_+\quad $$
and $$  D_{\bf i}^{-}(\ell_{-}(f))<\ell_-(h_-)\leq\ell_+(h_-)\leq\ell_+(f) \; \mbox{ for }  {\bf i} \in \mc I_-.$$
\end{lemma}
\begin{proof}
We establish the result with the positive indices (the case with negative indices is analogous). Write $f=\sum_{i=0}^{d} g_i (\ell_{+}(f))^i$ with $\ell_{+}(g_i)<\ell_{+}(f)$. Then, for any $M\in \MM_+$, using the fact that $\{-,M\}:S(\gg)\to S(\gg)$ is a derivation we get
$$D_M(f):=\{f,M\}=s_{+,f}\{\ell_+(f),M\}+\sum_{i=0}^{d} \{g_i,M\} (\ell_+(f))^ i.$$
If $D_{M}^+(\ell_+(f))=\ell_+(\{\ell_+(f), M\})\in\LL_{+}(\ell_+(f))$, then, by definition, 
$$\ell_+(\{g_i,M\})<D_{M}^+(\ell_+(f)).$$ 
Also, since the order $<$ is compatible with the grading we get $\ell_+(f)<D_{M}^+(\ell_+(f))$, which also implies that $\ell_{+}(s_{+,f})<D_{M}^+(\ell_+(f))$. Thus, the upper-leader of $\{f,M\}$ is $D_{M}^+(\ell_+(f))$, and 
$$D_M(f)=\alpha s_{+,f}D_{M}^{+}(\ell_{+}(f))+h_{+}$$
 where $\alpha$ is the coefficient of $D_{M}^+(\ell_{+}(f))$ when writing $\{\ell_+(f),M\}$ in terms of the basis $\MM$, and $h_+=\sum_{i=0}^{d} \{g_i,M\} (\ell_+(f))^ i-s_{+,f}\left(\{\ell_+(f),M\}-\alpha D_{M}^{+}(\ell_{+}(f))\right)$. 

Finally, note that all terms of $h_+$ are of the form $e$ or appear in $\{e,M\}$ for some $e\in \MM$ that appears in $f$. Since the order $<$ on $\mf M$ is compatible with the grading and $M \in \mf M_+$, we get that all these terms are larger or equal to $\ell_-(f)$ showing that $\ell_-(f)\leq \ell_-(h_+)$.

We have thus establish the result for the case when ${\bf i}$ is the 1-tuple $M$. For longer length tuples simply iterate this process.
\end{proof}

\begin{definition}\
\begin{enumerate}
\item [(i)] Let $f,g\in S(\gg)$ be nonconstant (i.e.,  not in $\kk$). We say that $g$ is partially reduced with respect to $f$ if no element in $\LL_{\pm}(\ell_{\pm}(f))$ appears (nontrivially) in $g$. If in addition $\ell_+(f)$ appears in $g$ only with degree $< d_{+,f}$, we say that $g$ is reduced with respect to $f$. 
\bigskip

\item [(ii)] Suppose $\Lambda=(f_i)_{i=1}^n$ is a sequence of nonconstant elements of $S(\gg)$ with $n\leq \omega$. If $g$ is in $S(\gg)$, we say that $g$ is (partially) reduced with respect to $\Lambda$ if $g$ is (partially) reduced with respect to every element in $\Lambda$  (by convention constant elements are reduced).
Furthermore, we say that $\Lambda$ is (partially) reduced if every $f_j$ is (partially) reduced with respect to $f_i$ for all $i<j$. 
\end{enumerate}
\end{definition}

\begin{remark}
Note that in a reduced sequence of elements from $S(\gg)$ distinct elements have distinct upper-rank (where recall that $rk_+(f)=(\ell_+(f),d_{+,f})$).
\end{remark}

\

The following Elimination Algorithm is one of the key ingredients of the proof of Theorem \ref{basisPoisson}.

\begin{theorem}\label{elimination}
Let $\Lambda$ be a finite sequence  of nonconstant elements of $S(\gg)$ and $g\in S(\gg)$. 
\begin{enumerate}
\item There exists $g_1\in S(\gg)$ partially reduced with respect to $\Lambda$ and integers $r_f\geq 0$ for each $f\in \Lambda$ such that
$$\prod_{f\in \Lambda}s_f^{r_f} g\;\equiv  \; g_1 \quad \operatorname{mod} [\Lambda]$$
and 
$$\min\{\ell_-(g),\ell_-(\Lambda)\}\leq \ell_{-}(g_1)\leq \ell_{+}(g_1)\leq \max\{\ell_+(g),\ell_+(\Lambda)\}.$$
\item Furthermore, if $\Lambda$ is reduced, then there also exists $g_0\in S(\gg)$ reduced with respect to $\Lambda$ and integers $m_f, n_f\geq 0$ for each $f\in \Lambda$ such that
$$\prod_{f\in \Lambda}i_{+,f}^{m_f}s_f^{n_f} g\;\equiv  \; g_0 \quad \operatorname{mod} [\Lambda]$$
and 
$$\ell_-(g_1)\leq \ell_-(g_0)\leq \ell_+(g_0)\leq \ell_+(g_1).$$
\end{enumerate}
\end{theorem}
\begin{proof}
(1) If no element in $\LL_+(\ell_+(f))$, for $f\in \Lambda$, appears in $g$ then we let $g_{+,1}=g$ and $r_{+,f}=0$. Otherwise, let $M$ be the largest element of $\mathfrak M$ that appears in $g$ and is in $\LL_+(\ell_+(f))$ for some $f\in \Lambda$. Then $M=D_{\bf i}^+(\ell_+(f))$ for some ${\bf i}=(M_1,\dots,M_n)\in \mathcal I_+$. By Lemma~\ref{separant}, we have
$$D_{\bf i}(f)=-\alpha s_{+,f}D_{\bf i}^+(\ell_+(f))+h_+$$
for some $\alpha\in \kk^*$ and $h_+\in S(\gg)$ with $\ell_-(f)\leq\ell_-(h_+)\leq \ell_+(h_+)<M$. 

Now write $g=\sum_{j=0}^r h_j M^j$ where $h_j$ is free of $M$ (and so the largest element of $\mathfrak M$ appearing in $h_j$ that is in $\LL_+(\ell_+(p))$ for some $p\in \Lambda$ is strictly less than $M$). Then, $g=\sum_{j=0}^r h_j (D_{\bf i}^+(\ell_+(f)))^j$ and so
\begin{equation}\label{eliminate}
s_{+,f}^rg=\sum_{j=0}^rh_js_{+,f}^{r-j}(s_{+,f} D_{\bf i}^+(\ell_+(f)))^j\equiv \sum_{j=0}^r h_js_{+,f}^{r-j}\, (\alpha^{-1}h_+)^j \quad \operatorname{mod}[\Lambda] .
\end{equation}
Since $\ell_+(s_{+,f})\leq \ell_+(f)<M$  and $\ell_+(h_+)<M$, the largest element from $\mathfrak M$ appearing in $g'_+ := \sum_{j=0}^r h_js_f^{r-j}\, (\alpha^{-1}h_+)^j$ that is in $\LL_+(\ell_+(p))$ for some $p\in \Lambda$ is strictly less than $M \leq \ell_+(g)$. 
Furthermore, since $\ell_-(f)\leq\ell_-(h_+)$, we have 
$$\min\{\ell_-(g),\ell_-(f)\} \leq \ell_-(g'_+)\leq \ell_+(g'_+)\leq \ell_+(g)\leq \max\{\ell_+(g),\ell_+(f)\}.$$

Repeat the above process on $g'_+$ until we reach $g_{+,1}\in S(\gg)$ such that no element in $\LL_+(\ell_+(f))$, for $f\in \Lambda$, appears in it (this process eventually terminates as $\Lambda$ is finite). Note that, by \eqref{eliminate}, 
$$\prod_{f\in \Lambda}s_{+,f}^{r_{+,f}} g\;\equiv  \; g_{+,1} \quad \operatorname{mod} [\Lambda],$$
for some $r_{+,f}$, and 
$$\min\{\ell_-(g),\ell_-(\Lambda )\} \leq \ell_-(g_{+,1}) \leq \ell_{+}(g_{+,1})\leq \max\{\ell_+(g),\ell_+(\Lambda)\}.$$

Now, if no element in $\LL_-(\ell_-(f))$, for $f\in \Lambda$, appears in $g_{+,1}$ then we let $g_1=g_{+,1}$ and $r_f=r_{+,f}$. Otherwise, we perform the counterpart (i.e., negative-indices) of the above process to $g_{+,1}$. Again the process will eventually terminate, as $\Lambda$ is finite, and yields the desired $g_1$.

\medskip 

(2) Let $\Lambda=(f_1,\dots,f_s)$ and now we assume that it is reduced. Since distinct elements in $\Lambda$ have distinct upper-rank, there are $k_1,\dots,k_s\in\{1,\dots,s\}$ such that 
$$rk_+(f_{k_1})>rk_+(f_{k_2})>\cdots >rk_+(f_{k_s}).$$
If $\ell_+(f_{k_1})$ does not appear in $g_1$ then we let $g_{0,1}=g_1$. Otherwise, assume $\ell_+(f_{k_1})$ appears in $g_1$ with degree $r\geq d:=d_{+,f_{k_1}}$. Write $g_1=\sum_{j=0}^r h_j (\ell_{+}(f_{k_1}))^j$ with $h_j$ free of $\ell_+(f_{k_1})$. Then, in 
\begin{equation}\label{initial}
i_{+,f_{k_1}}g_1-h_r(\ell_+(f_{k_1}))^{r-d}f_{k_1}
\end{equation}
$\ell_+(f_{k_1})$ appears with degree $\leq r-1$. Repeating this process yields $\tilde g_{0,1}$ reduced with respect to $f_{k_1}$. Since $\Lambda$ is reduced and $g_1$ is partially reduced with respect to $\Lambda_{k_1}:=(f_1,\dots,f_{k_1})$, $\tilde g_{0,1}$ is partially reduced with respect to $\Lambda_{k_1}$. However, the above process \eqref{initial} might yield $\tilde g_{0,1}$ that is not partially reduced with respect to $\Lambda_{k_1}^*:=\{f_{k_1+1},\dots,f_s\}$ (as $f_{k_1}$ is not necessarily partially reduced with respect to $\Lambda_{k_1}^*$). Note that since $g_1$ is partially reduced with respect to $\Lambda$, the only way that the above process could produce an element which is  not partially reduced with respect to $\Lambda_{k_1}^*$ is if in $f_{k_1}$ appears an element from $\LL_+(\ell_+{f_i})$ for some $i>k_1$. In fact, any element from $\LL_+(\ell_+(f_i))$, for some $i>k_1$, that appears in $\tilde g_{0,1}$ is $<\ell_+(f_{k_1})$ as $\ell_+(f_{k_1})$ cannot be in any $\ell_+(f_i)$ since it appears in $g_1$.

Now perform the algorithm from part (1) to $\tilde g_{0,1}$ with $f = f_{k_1}$. Notice that on the right-handed term in \eqref{eliminate} the upper-leaders of $s_{+,f}$ and $h_+$ will be $<\ell_{+}(f_{k_1})$ and in the $h_j$'s this basis term will appear with degree $<d_{+,f_{k_1}}$. Thus, the output of the algorithm from part (1) is $g_{0,1}$ with degree in $\ell_+(f_{k_1})$ strictly less than $d_{+,f_{k_1}}$. In other words, $g_{0,1}$ is partially reduced with respect to $\Lambda$ and reduced with respect to $f_{k_1}$.

Perform the same process with $f_{k_2}$ and $g_{0,1}$ to obtain $g_{0,2}$ partially reduced with respect to $\Lambda$ and reduced with respect to $f_{k_2}$. Note that this $g_{0,2}$ will also be reduced with respect to $f_{k_1}$. Indeed, if $\ell_+(f_{k_1})=\ell_+(f_{k_2})$ then it is clear since $rk_+(f_{k_1})>rk_+(f_{k_1})$; otherwise, $\ell_+(f_1)>\ell_+(f_{k_2})$ and in this case the degree of $\ell_+(f_{k_1})$ is not increased in the algorithm of part (1) applied to $g_{0,1}$, meaning that its degree in $g_{0,2}$ is $<d_{+,f_{k_1}}$ and so $g_{0,2}$ is reduced with respect to $f_{k_1}$. Repeating the above process yields the desired $g_0$.
 \end{proof}

The Elimination Algorithm yields the following useful corollary:

\begin{corollary}\label{cor:fingen}
Suppose that the basis $\MM$ of $\mf g$ satisfies condition \eqref{cofinite} on p. \pageref{cofinite}. 
Then, for any nonzero prime Poisson ideal $P$ of $S(\mf g)$, there is $h \in S(\mf g) \setminus P$ such that the localisation $(S(\mf g)/P)_h$ is an affine $\kk$-algebra.
\end{corollary}
\begin{proof}
Condition \eqref{cofinite} is equivalent to the hypothesis that for each $M\leq N\in \MM$ the $\kk$-subspace of $\mf g$ spanned by $\LL_-(M)\cup \LL_+(N)$ has finite codimension. 
Let $f$ be a nonzero element of $P$ such that $s_f$ is not in $P$ (for instance, choose $f\in P$ of minimal total degree). By part (1) of the Elimination Algorithm, for any $g\in S(\mf g)$ there is $g_1$ partially reduced with respect to $f$ and an integer $r\geq 0$ such that
$$s_f^r g\; \equiv \;  g_1 \quad \mod [f].$$
By our assumption, $g_1$ lives in the affine $\kk$-algebra generated by the finite set $\MM\setminus (\LL_-(\ell_-(f))\cup \LL_+(\ell_+(f)))$. Thus, setting $h=s_f$, we get that $(S(\mf g)/P)_h$ is a finitely generated $\kk$-algebra.
\end{proof}

Under the hypothesis of Corollary~\ref{cor:fingen} $\mf g$ is automatically Dicksonian by Remark~\ref{rem:newcond}.  Note that $S(\Vir)/(z) \cong S(W)$ is not contained in any affine 
{
(commutative)} $\kk$-algebra, so for the proof of  Corollary~\ref{cor:fingen} we need \eqref{cofinite} to hold for all $M \in \MM$.

\medskip

The following lemma is used to relate the Dicksonian condition to chain conditions in symmetric algebras.  

\begin{lemma}\label{finitereduced}
If $\gg$ is Dicksonian,
then every reduced sequence of $S(\gg)$ is finite. 
\end{lemma}
\begin{proof}
Suppose there is an infinite reduced sequence $\Lambda=(f_1,f_2,\cdots)$ in $S(\gg)$. We claim that an infinite subsequence of $(\ell_-(f_i),\ell_+(f_i))_{i=1}^\infty$ is leading-Dicksonian. This sequence of pairs clearly satisfies 
$\ell_-(f_i)\leq \ell_+(f_i)$. Also, by definition of reduced sequence, it satisfies
$$\ell_-(f_j)\notin\LL_-(\ell_-(f_i)) \quad \text{ and } \quad \ell_+(f_j)\notin\LL_+(\ell_+(f_i))$$
for all $i<j$. Thus, the only condition from the definition of leading-Dicksonian that is missing is that the elements of $(\ell_-(f_i),\ell_+(f_i))_{i=1}^\infty$ are distinct. Now, because in the definition of reduced sequence we require that $\ell_+(f_i)$ appears in $f_j$ only with degree strictly less than $d_{+,f}$, equality of $\ell_+(f_j)$ and $\ell_+(f_i)$ with $i<j$ can only happen finitely many times. Thus there will be a subsequence with the desired properties, contradicting our assumption. 
\end{proof}

Recall that if $I$ is a radical ideal then the division ideal $I:s$ is again radical, and if $I$ is radical and Poisson then the same is true of $I:s$.
Further, if $I$ is Poisson then the ideal $$I:s^\infty=\{f\in S(\gg): fs^n \in I \text{ for some } n \geq 0\}$$ is again Poisson, by a similar argument to the proof that $\CP$ is  a radical conservative system (see \S\ref{specialPoi}).

\medskip

For any finite subset $\Lambda$ of $S(\gg)$, we let $[\Lambda]$ denote the Poisson ideal generated by $\Lambda$ in $S(\gg)$ and set $i_+s_\Lambda=\prod_{f\in \Lambda}i_{+,f}s_f$, where recall that $s_f=s_{+,f}\cdot s_{-,f}$. In Theorem~\ref{finiteprime} below we will be looking at Poisson ideals of the form $[\Lambda]:i_+s_\Lambda^\infty$.

\bigskip

The above Elimination Algorithm yields the following fundamental ``basis theorem'' for prime Poisson ideals:

\begin{theorem}\label{finiteprime}
Assume $\gg$ is Dicksonian.
 If $P$ is a nonzero prime Poisson ideal of $S(\gg)$, then there is an reduced set $\Lambda$ (which is hence finite by Lemma~\ref{finitereduced}) contained in $P$ such that $i_+s_\Lambda\notin P$ and
$$P=[\Lambda]:i_+s_\Lambda^\infty .$$
\end{theorem}

\begin{proof}
Let $f_1$ be a nonzero element of $P$ of minimal total degree. Since the upper-initial $i_{+,f_1}$, the upper-separant $s_{+,f_1}$ and the lower-separant $s_{-,f_1}$ are nonzero of total degree smaller than that of $f_1$,  none of them is in $P$. Let $\Lambda_1$ be the singleton sequence $(f_1)$. As $P$ is prime it follows that $i_+s_{\Lambda_1}\notin P$, and so
$$[\Lambda_1]:i_{+}s_{\Lambda_1}^\infty \subseteq P.$$ 
If $P=[\Lambda_1]:i_{+}s_{\Lambda_1}^\infty$ we are done. Otherwise, there is $g\in P$ but not in $[\Lambda_1]:i_{+}s_{\Lambda_1}^\infty$. By Theorem \ref{elimination}, there is $g_0$ which is reduced with respect to $\Lambda_1$ and integers $m_{f_1},n_{f_1}\geq 0$ such that
$$i_{+,f_1}^{m_{f_1}}s_{f_1}^{n_{f_1}} g\;\equiv  \; g_0 \quad \operatorname{mod} [\Lambda_1 ].$$
So $g_0$ is in $P$ and is nonzero (otherwise $g$ would be in $[\Lambda_1]:i_+s_{\Lambda_1}^\infty$). So we can choose $f_2$ of minimal total degree among the nonzero elements in $P$ that are reduced with respect to $\Lambda_1$. Then, the upper initial and upper and lower separants of $f_2$ are not in $P$ (as they are all nonzero reduced with respect to $\Lambda_1$ and of total degree smaller than that of $f_2$). Let $\Lambda_2$ be the sequence $(f_1,f_2)$. Then $\Lambda_2$ is reduced and, as $P$ is prime, $i_+s_{\Lambda_2}\notin P$.  So we have 
$$[\Lambda_2]:i_+s_{\Lambda_2}^\infty \subseteq P.$$
If $[\Lambda_2]:i_+s_{\Lambda_2}^\infty$ is not equal to $P$, we can continue this process and find $f_3$ such that the sequence $\Lambda_3=(f_1,f_2,f_3)$ is reduced and $i_+s_{\Lambda_3}\notin P$ . This process must eventually stop since reduced sequences are finite (by Lemma~\ref{finitereduced}). Thus, this process yields a (finite) reduced sequence $\Lambda$ such that $i_{+}s_{\Lambda}$ is not in $P$ and
$$P=[\Lambda]:i_+s_{\Lambda}^\infty.$$
\end{proof}

\subsection{Proof of main theorem}

We can now easily prove Theorem~\ref{basisPoisson}. We restate it for the reader's convenience. 

\medskip

\begin{theorem}\label{thm:basisPoissonrestate}
Let $\kk$ be a field of characteristic zero and $\gg$ a graded Lie $\kk$-algebra. If $\gg$ is Dicksonian, 
then the Poisson algebra $S(\gg)$ has ACC on radical Poisson ideals.
\end{theorem}

\begin{proof}
By Theorem~\ref{finiteprime}, for every prime Poisson ideal $P$ of $S(\gg)$ there is a finite reduced sequence $\Lambda\subset P$ such that $i_+s_{\Lambda}\notin P$ and 
\begin{equation}\label{primequal}
P=[\Lambda]:i_+s_\Lambda^\infty .
\end{equation}
Setting $s=i_+s_\Lambda$ we see that
$$[\Lambda]:i_+s_\Lambda^\infty\subseteq \{\Lambda\}:s$$
where $\{\Lambda\}$ denotes the radical Poisson ideal generated by $\Lambda$ in $S(\gg)$. Since $P$ is prime Poisson and $s\notin P$, we get $\{\Lambda\}:s\subseteq P$. By \eqref{primequal}, we actually have
$$P=\{\Lambda\}:s .$$
The result now follows immediately from Theorem~\ref{conservativePoisson}(i).
\end{proof} 

We finish this section by noting that our Poisson basis theorem applies to a wide collection of Poisson algebras.

\begin{corollary}
Let $A$ be a Poisson algebra over a field $\kk$ of characteristic zero and $N\subseteq A$ a Poisson-generating set (i.e., $N$ generates $A$ as a Poisson algebra). If the Lie $\kk$-algebra generated by $N$ 
is Dicksonian (with respect to some $\ZZ$-grading and some ordered basis of homogeneous elements) 
then $A$ has ACC on radical Poisson ideals. 
\end{corollary}
\begin{proof}
Theorem~\ref{basisPoisson} establishes the ACC (on radical Poisson ideals) for the Poisson algebra $S(\mathfrak g)$ where $\gg$ is the Lie algebra generated by $N$. But $A$ is a factor of the symmetric algebra $S(\mathfrak g)$ by a Poisson ideal, and hence we also have the ACC on radical Poisson ideals for $A$. 
\end{proof}

\section{Examples}\label{examples}

In this section we provide a wide family of examples of graded Lie algebras 
which are Dicksnonian. 
Throughout this section $\kk$ is a field of characteristic zero. At the end, we prove that all simple graded Lie algebras of polynomial growth 
are Dicksonian.
Consequently, by Theorem~\ref{basisPoisson},  their symmetric algebras have ACC on radical Poisson ideals (namely, we prove Corollary \ref{simplecase}) and, by Theorem~\ref{conservativePoisson}(ii), every radical Poisson ideal is a finite intersection of prime Poisson ideals. 

\bigskip

\subsection{The Witt algebra}

In Section \ref{leading} we saw that the Witt algebra $W$, the positive Witt algebra $W_+$,  the Cartan algebra $\WW_1$, and the Virasoro algebra $\Vir$
are Dicksonian 
(with respect to the natural choice of ordered basis). 

We thus get the following consequences from Theorems \ref{basisPoisson} and  \ref{conservativePoisson}(ii).

\begin{corollary}
Let $\gg$ be one of $W$, $W_+$,  $\WW_1$, or $\Vir$. Then $S(\gg)$ has ACC on radical Poisson ideals and every radical Poisson ideal is a finite intersection of prime Poisson ideals. 
\end{corollary}

\bigskip

\subsection{Cartan algebras}\label{Cartan}

Let $n\geq 2$. In this section we prove that the Cartan algebra $\WW_n$ has a basis with an order $<$ such that $\WW_n$ has no infinite leading-Dicksonian sequence. Theorem~\ref{basisPoisson} then tells us that $S(\WW_n)$ has ACC on radical Poisson ideals. 

\medskip

Let $\kk$ be a field of characteristic zero.  
Recall that the {\em Cartan algebra} $\WW_n$ over $\kk$  is the Lie $\kk$-algebra of derivations on $\kk[x_1,\cdots,x_n]$. We will use multi-index notation; that is, for $i=(i_1,\dots,i_n)\in \mathbb N^n$ we write 
$$x^i=(x^{i_1},\dots,x^{i_n}).$$
By $|i|$ we mean $i_1+\cdots+i_n$. Also, when we write $i\leq j$ with $i,j\in \mathbb N^n$ we mean $i_1\leq j_1,\dots,i_n\leq j_n$ (i.e., $i< j$ means $i$ is less than $j$ in the product order of $\mathbb N^n$). For $k\in \{1,\dots,n\}$, we let $1_k$ denote the $n$-tuple with a $1$ in the $k$-entry and zeroes elsewhere. 

We choose as $\MM$ the natural basis for $\WW_n$; that is, 
$$\{x^i\partial_k: i\in \mathbb N^n \text{ and } k\in \{1,\dots,n\}\}. $$
where $\partial_k=\frac{\partial}{\partial x_k}$. We equip this basis with the following (total) ordering: 
$$x^i\partial_k < x^j\partial_\ell \quad \Leftrightarrow\quad (|i|,k,i_n,\dots,i_1)<_{\text{lex}}(|j|,\ell,j_n,\dots,j_1) .$$
So this ordering is compatible with the natural $\ZZ$-grading of $\WW_n$. 

We now check that, with respect to this order, $\WW_n$ 
is Dicksonian.
First we need a lemma. A word on notation.  In the lemma below we write $0 = (0,0, \dots, 0)$, and if  $k=1$ we set $(i_1,\dots,i_{k-1})=0$ for convenience of exposition.

\begin{lemma}\label{leadingW}
Let $i\in \mathbb N^n$ and $k\in \{1,\dots,n\}$. 
\begin{enumerate}
\item [(i)] If $(i_1,\dots,i_{k-1})\neq 0$, then $\LL_+(x^i\partial_k)$ contains
$$\{x^r\partial_k :\; r>i\} .$$
\item [(ii)] If $(i_1,\dots,i_{k-1})=0$, then $\LL_+(x^i\partial_k)$ contains
$$\{x^r\partial_k :\; r>i, (r_1,\dots,r_{k-1})=0, \text{ and } r_k\neq 2i_k-1\} .$$
If in addition $i_k=1$ and $i\neq 1_k$, then $\LL_+(x^i\partial_k)$ contains
$$\{x^r\partial_k :\; r>i, (r_1,\dots,r_{k-1})=0\}$$
\end{enumerate}
\end{lemma}

\begin{proof}
(i) Assume $(i_1,\dots,i_{k-1})\neq 0$, say $i_m\neq 0$ with $m\in\{1,\dots,k-1\}$. Then, for any $j\in \mathbb N^n$,
$$[x^i\partial_k,x^j\partial_m]=j_kx^{i+j-1_k}\partial_m - i_mx^{i+j-1_m}\partial_k.$$
So, since $i_m\neq 0$, for any $r>i$ if we set $j=r-i+1_m$ we get
$$\ell_+([x^i\partial_k,x^j\partial_m])=x^r \partial_k.$$
Moreover, if $x^u\partial_\ell<x^i\partial_k$, one can easily  see that 
$$\ell_+([x^u\partial_\ell,x^j\partial_m])<\ell_+([x^i\partial_k,x^j\partial_m]).$$
Thus, $\LL_+(x^i\partial_k)$ contains all elements of the form $x^r\partial_k$ for $r>i$.

\

(ii) Assume $(i_1,\dots,i_{k-1})=0$. Let $j\in\mathbb N^n$ such that $(j_1,\dots,j_{k-1})=0$. Then 
\begin{equation}\label{firsteq}
[x^i\partial_k, x^j\partial_k]=(j_k-i_k)x^{i+j-1_k}\partial_k
\end{equation}
while for any $\ell<k$ and $u\in \mathbb N^n$
\begin{equation}\label{secondeq}
[x^u\partial_\ell,x^j\partial_k]=-u_kx^{u+j-1_k}\partial_\ell.
\end{equation}
For any $r>i$ with $(r_1,\dots,r_{k-1})=0$ and  $r_k\neq 2i_k-1$, if we set $j=r-i+1_k$ then by \eqref{firsteq} we get
$$\ell_+([x^i\partial_k, x^j\partial_k])=x^{r}\partial_k$$
and by \eqref{secondeq}, for any $x^u\partial_\ell<x^i\partial_k$, we see that 
$$\ell_+([x^u\partial_\ell,x^j\partial_k])<\ell_+([x^i\partial_k,x^j\partial_k]).$$
Thus, $\LL_+(x^i\partial_k)$ contains all elements of the desired form.

Now, if in addition $i_k=1$ and $i\neq 1_k$, we must show that $\LL_+(x^i\partial_k)$ also contains all elements of the form $x^r\partial_k$ where $r>i$, $(r_1,\dots,r_{k-1})=0$ and $r_k=1$. Note that since $i\neq 1_k$, we must have that $k<n$ and there is $m>k$ such that $i_m\neq 0$. Let $j=r-i+1_{m}$, which is in $\NN^k$ since $r>i$ in the product order.
Then
$$\ell_+([x^i\partial_k,x^j\partial_{m}])=x^{r}\partial_k.$$
Moreover, for any $x^u\partial_\ell<x^i\partial_k$, we see that 
$$\ell_+([x^u\partial_\ell,x^j\partial_m])<\ell_+([x^i\partial_k,x^j\partial_m]).$$
Thus, $\LL_+(x^i\partial_k)$ contains such $x^r\partial_k$.
\end{proof}

We now recall Dickson's lemma, as we will make use of it. Let $\mathcal N=\mathbb N^n\times\{1,\dots,n\}$. We equip $\mathcal N$ with the following ordering: for $(i,k),(j,\ell)\in \mathcal N$, we set $(i,k)<(j,\ell)$ if and only if $k=\ell$ and $i<j$ (the latter denotes the product order of $\mathbb N^n$). A sequence of elements $(a_i)$ from $\mathcal N$ is called \emph{Dicksonian} if $a_j\not\geq a_i$ for $j>i$.

\begin{theorem}[Dickson's lemma \cite{Fig}]\label{Dickson}
Every Dicksonian sequence of $\mathcal N$ is finite.
\end{theorem}

\begin{proposition}\label{forW}
With respect to the above ordering on $\MM$, the Cartan algebra $\WW_n$
is Dicksonian.
\end{proposition}
\begin{proof}
Assume towards a contradiction that there is an infinite leading-Dicksonian sequence $((M_i,N_i))_{i=1}^\infty$ in $\WW_n$. From the sequence of $N_i$'s we can find an infinite sequence $(a_i)$ of elements in $\MM_+$ such that $a_j\notin \LL_+(a_i)$ for $i<j$. To each $a_i=x^j\partial_k\in \MM_+$ we associate the element $b_i:=(j,k)\in \mathcal N$. This gives us an infinite sequence $(b_i)$ of $\mathcal N$.  By Lemma~\ref{leadingW}, there is an infinite subsequence of $(b_i)$ which is Dicksonian, but this contradicts Dickson's lemma, Theorem~\ref{Dickson}.
\end{proof}

The proof of Proposition~\ref{forW} is the origin of our use of the term {\em Dicksonian} to describe our key condition on Lie algebras.

\medskip

\subsection{Special Cartan algebras}\label{specialCartan}

Let $n\geq 2$. Recall that the \emph{special Cartan algebra} $\SS_n$ is the Lie subalgebra of $\WW_n$ given by elements of the form
\begin{equation}\label{SSn-def}
p_1\partial_1+\cdots + p_n\partial_n \quad \text{ such that }\quad \partial_1(p_1)+\cdots+\partial_n(p_n)=0
\end{equation}
where $p_i\in \kk[x_1,\dots,x_n]$.

In this section we prove that $\SS_n$ is Dicksonian and hence, by Theorem~\ref{basisPoisson}, the symmetric algebra $S(\SS_n)$ has ACC on radical Poisson ideals. 

Induce the grading on $\SS_n$ from the natural grading on $\WW_n$.  We seek an ordered homogeneous basis  of $\SS_n$.  Let $\mf N$ be the subset of $\SS_n$ consisting of elements of the form
$$x^i\partial_1, \quad \text{ such that }  i\in \mathbb N^n \text{ and } i_1=0,$$
together with the elements of the form
$$i_kx^{i-1_k}\partial_1 - i_1x^{i-1_1}\partial_k, \quad \text{ with }  2\leq k\leq n, i\in \mathbb N^n \text{ and } i_1\neq 0.$$

\begin{lemma}
The set $\mf N$ is a $\kk$-basis for $\SS_n$.
\end{lemma}
\begin{proof}
A straightforward computation shows that $\mf N$ is $\kk$-linearly independent and contained in $\SS_n$. To see that $\mf N$ $\kk$-spans $\SS_n$, 
note that $\del_1$ defines a surjective 
{
linear} endomorphism of $\kk[x_1, \dots, x_n]$ whose kernel is $\kk[x_2, \dots, x_n]$.
Thus $\del_1$ is split by the map
$\del_1^{-1} :\kk[x_1, \dots, x_n]  \to x_1 \kk[x_1, \dots, x_n] $ defined by extending
$ x^i \mapsto \frac{x_1 x^i}{i_1+1}$ linearly.
  Solutions to \eqref{SSn-def} are thus given by
\[ p_1 \in \del_1^{-1}(-\del_2 (p_2) -\dots - \del_n(p_n)) + \kk[x_2, \dots, x_n],\]
and all such solutions are clearly in the $\kk$-span of $\mf N$. 
\end{proof}

We now equip $\mf N$ with the  (total) order $<$ induced from the order we defined on our basis $\MM$ of $\WW_+$. That is, let $\ell_{+, \MM}$ denote the leading term of an element of $\mf N$ with respect to $\MM$ and define $M <_{\mf N} N$ if and only if $\ell_{+, \MM}(M) <_{\MM} \ell_{+, \MM}(N)$.  
Explicitly, let
 $i,j\in \mathbb N^n$ and $k,l\in \{2,\dots,n\}$. When $i_1=j_1=0$ we set
$$x^i\partial_1< x^j\partial_1  \quad \Leftrightarrow\quad (|i|,i_n,\dots,i_2)<_{\text{lex}}(|j|,j_n,\dots,j_2),$$
when $i_1=0$ and $j_1\neq 0$
$$x^i\partial_1< j_k x^{j-1_k}\partial_1 - j_1 x^{j-1_1}\partial_k  \quad \Leftrightarrow\quad (|i|,1)<_{\text{lex}}(|j|-1,k),$$
and when $i_1\neq 0$ and $j_1\neq 0$
\begin{multline*}
i_kx^{i-1_k}\partial_1- i_1 x^{i-1_1}\partial_k <  j_l x^{j-1_l}\partial_1 - j_1 x^{j-1_1}\partial_l  \quad \Leftrightarrow \\
 (|i|,k,i_n,\dots,i_1)<_{\text{lex}}(|j|,l,j_n,\dots,j_1).
 \end{multline*}

\medskip

We now prove a lemma which can be thought of as the $\SS_n$ analogue of  Lemma~\ref{leadingW}.

\begin{lemma}\label{leadingS}
Let $i\in \mathbb N^n$ and $k\in\{2,\dots,n\}$.
\begin{enumerate}
\item [(i)] If $i_1= 0$, then $\LL_+(x^i\partial_1)$ contains
$$\{x^r\partial_1 :\; r>i \text{ and } r_1=0\} .$$
\item [(ii)] If $i_1=1$, then $\LL_+(i_kx^{i-1_k}\partial_1 - i_1x^{i-1_1}\partial_k)$ contains
$$\{r_kx^{r-1_k}\partial_1 - r_1x^{r-1_1}\partial_k:\; r>i \text{ and } \; \sum_{j=2}^n(r_j-i_j)\neq 1 \} .$$ \medskip
If $i_1\geq 2$, then $\LL_+(i_kx^{i-1_k}\partial_1 - i_1x^{i-1_1}\partial_k)$ contains
$$\{r_kx^{r-1_k}\partial_1 - r_1x^{r-1_1}\partial_k:\; r>i, \; \sum_{j=2}^n(r_j-i_j)\neq 1 \text{ , and }\;  r_k(r_1-i_1+2)+1\neq i_1\} .$$
\end{enumerate}
\end{lemma}

\begin{proof}
(i) Assume $i_1=0$ and $r>i$ with $r_1=0$. If $j=r-i+1_1+1_n$, then 
$$[x^i\partial_1,j_nx^{j-1_n}\partial_1 - j_1x^{j-1_1}\partial_n]=(i_n+j_n) x^{r}\partial_1 .$$
Since $j_n>0$, it follows that
$$\ell_+([x^i\partial_1,j_nx^{j-1_n}\partial_1 - j_1x^{j-1_1}\partial_n])=x^r\partial_1 .$$
Moreover, if $x^u\partial_1<x^i\partial_1$, one easily  checks that 
$$\ell_+([x^u\partial_1,j_nx^{j-1_n}\partial_1 - j_1x^{j-1_1}\partial_n])<\ell_+([x^i\partial_1,j_nx^{j-1_n}\partial_1 - j_1x^{j-1_1}\partial_n]),$$
and if $u_kx^{u-1_k}\partial_1-u_1x^{u-1_1}\partial_k < x^i\partial_1$ one also easily checks that
\begin{multline*}
\ell_+([u_kx^{u-1_k}\partial_1-u_1x^{u-1_1}\partial_k, j_nx^{j-1_n}\partial_1 - j_1x^{j-1_1}\partial_n]) \\
< \ell_+([x^i\partial_1,j_nx^{j-1_n}\partial_1 - j_1x^{j-1_1}\partial_n]).
 \end{multline*}
Thus, $\LL_+(x^i\partial_1)$ contains all elements of the form $x^r\partial_1$ for $r>i$ with $r_1=0$.

\bigskip

(ii) Let $i\in \mathbb N^n$ and $r>i$. 
\smallskip

Set $r'=(r_1,i_2,\dots,i_n)$. If we let $j=(r_1-i_1+1)\cdot 1_1+1_k$, then
$$[i_kx^{i-1_k}\partial_1-i_1x^{i-1_1}\partial_k, j_kx^{j-1_k}\partial_1-j_1x^{j-1_1}\partial_k]$$
equals
\begin{equation}\label{here}
(i_k(r_1-i_1+1)-i_1)\cdot\left(r'_k x^{r'-1_k}\partial_1- r'_1 x^{r'-1_1}\partial_k\right) .
\end{equation}
On the other hand, if $i_1\geq 1$ set $r''=(i_1-1, r_2,\dots,r_n)$. If we let $j=(0,r_2-i_2,\dots,r_n-i_n)$, we get
\begin{equation}\label{there}
[i_kx^{i-1_k}\partial_1-i_1x^{i-1_1}\partial_k, x^j\partial_1]=(-i_1)\left( r''_kx^{r''-1_k}\partial_1 - r''_1x^{r''-1_1}\partial_k \right) .
\end{equation}
Note that $x^j\partial_1\in \mf N_+$ as long as $\sum_{j=2}^n(r_j-i_j)> 1$.
\medskip

Now, in case $i_1=1$, from \eqref{here} and since $2i_k\neq 1$ we see that 
$$\ell_+([i_kx^{i-1_k}\partial_1-i_1x^{i-1_1}\partial_k,x^{2\cdot 1_1}\partial_1-2x^{1_1+1_k}\partial_k])=i_kx^{i+1_1-1_k}\partial_1-(i_1+1)x^{i}\partial_k$$
and if $N<i_kx^{i-1_k}\partial_1-i_1x^{i-1_1}\partial_k$ then
$$\ell_+([N,x^{2\cdot 1_1}\partial_1-2x^{1_1+1_k}\partial_k])<\ell_+([i_kx^{i-1_k}\partial_1-i_1x^{i-1_1}\partial_k,x^{2\cdot 1_1}\partial_1-2x^{1_1+1_k}\partial_k]).$$
So 
$$i_kx^{i+1_1-1_k}\partial_1-(i_1+1)x^{i}\partial_k$$
is in $\LL_+(i_kx^{i-1_k}\partial_1-i_1x^{i-1_1}\partial_k)$. By a similar argument, using now \eqref{there} one checks that
$$r_kx^{(1,r_2,\dots,r_k-1,\dots,r_n)}\partial_1-x^{(0,r_2,\dots,r_n)}\partial_k$$
is in $\LL_+(i_kx^{i-1_k}\partial_1-i_1x^{i-1_1}\partial_k)$. Finally, another application of \eqref{here} yields that
$$r_kx^{r-1_k}\partial_1 - r_1x^{r-1_1}\partial_k$$
is in $\LL_+(i_kx^{i-1_k}\partial_1-i_1x^{i-1_1}\partial_k)$, as desired.

\medskip

Now, in the case $i\geq 2$ the argument is similar to the one above. Using \eqref{there} one checks that
$$r_kx^{(i_1-1,r_2,\dots,r_k-1,\dots,r_n)}\partial_1-(i_1-1)x^{(i_1-2,r_2,\dots,r_n)}\partial_k$$
is in $\LL_+(i_kx^{i-1_k}\partial_1-i_1x^{i-1_1}\partial_k)$. And finally, using \eqref{here} one checks that
$$r_kx^{r-1_k}\partial_1 - r_1x^{r-1_1}\partial_k$$
is in $\LL_+(i_kx^{i-1_k}\partial_1-i_1x^{i-1_1}\partial_k)$ as long as $r_k(r_1-i_1+2)+1\neq i_1$, as desired.

\end{proof}

\begin{proposition}\label{forS}
With respect to the above ordering on $\mf N$, the special Cartan algebra $\SS_n$ 
is Dicksonian.
\end{proposition}
\begin{proof}
The proof is almost identical to the proof of Proposition~\ref{forW} (but using Lemma~\ref{leadingS} rather than Lemma~\ref{leadingW}). Namely, to a given infinite leading-Dicksonian sequence of $\SS_n$ one naturally associates an infinite sequence $(b_i)$ of $\mathcal N$ (recall from \S\ref{Cartan} that the latter denotes $\mathbb N^n\times\{1,\dots,n\}$ equipped with the natural product order). This is done as follows: for $i\in \mathbb N^n$ and $k\in \{2,\dots,n\}$; when $i_1=0$
$$x^i\partial_1 \quad \mapsto \quad (i,1)$$
and when $i_1\neq 0$
$$i_kx^{i-1_k}\partial_1 -i_1x^{i-1_1}\partial_k \quad \mapsto \quad (i-1_1,k) .$$
Now by Lemma~\ref{leadingS} there is an infinite subsequence of $(b_i)$ which is Dicksonian, contradicting Dickson's lemma, Theorem~\ref{Dickson}. 
\end{proof}

\medskip

\subsection{Hamiltonian Cartan algebras}
Let $n=2m$ with $m$ a positive integer. Recall that the \emph{Hamiltonian Cartan algebra} $\HH_n$ is the Lie subalgebra of $W_n$ given by elements of the form
$$D_H(p):=\sum_{\ell=1}^{m}\left( \partial_{m+\ell}(p)\partial_\ell \; -\; \partial_\ell(p)\partial_{m+\ell} \right)$$
for $p\in \kk[x_1,\dots,x_n]$. In fact $D_H:\kk[x_1,\dots,x_n]\to W_n$ as defined above is a  $\kk$-linear mapping with kernel $\kk$. 
One can easily derive that for $p,q\in \kk[x_1,\dots,x_n]$ we have $[D_H(p),D_H(q)]=D_H(h)$ where 
\begin{equation}\label{forhamilton}
h=\sum_{\ell=1}^{m}\left(\partial_{m+\ell}(p)\partial_\ell(q) \; -\; \partial_\ell(p)\partial_{m+\ell}(q)\right)
= D_H(p)(q).
\end{equation}
In other words, $[D_H(p),D_H(q)]=D_H(\{p,q\})$ where $\{p,q\}=D_H(p)(q)$.

In this section we prove that $\HH_n$ 
is Dicksonian.
Hence, by Theorem~\ref{basisPoisson}, the symmetric algebra $S(\HH_n)$ has ACC on radical Poisson ideals. 

\smallskip

Again, we seek an ordered homogeneous basis $\MM$ for $\HH_n$.  
Let $\MM=\{D_H(x^i): i\in \mathbb N^n \text{ and } i\neq {\bf 0}\}\subset \HH_n$. The set $\MM$ is clearly a $\kk$-basis for $\HH_n$. From~\eqref{forhamilton}, we see that
\begin{equation}\label{basishamilton}
[D_H(x^i),D_H(x^j)]=D_H\left(\sum_{\ell=1}^m \left( i_{m+\ell} j_\ell - i_\ell j_{m+\ell} \right)x^{i+j-1_\ell -1_{m+\ell}}\right) .
\end{equation}

We equip $\MM$ with the following (total) order. Given nonzero $i,j\in \mathbb N^n$, we set 
$$D_H(x^i)< D_H(x^j) \quad \Leftrightarrow\quad (|i|,i_n,\dots,i_1)<_{\text{lex}}(|j|,j_n,\dots,j_1) . $$

We now prove a lemma which can be thought of as the $\HH_{n}$ analogue of  Lemmas~\ref{leadingW} and \ref{leadingS}.

\begin{lemma}\label{leadingH}
Let $i\in \mathbb N^n$ and $1\leq \ell \leq m$. 
\begin{enumerate}
\item If $i_\ell \neq 2 \,i_{m+\ell}$, then $\mathcal L_+(D_H(x^i))$ contains
$$\{D_H(x^{i+r\cdot 1_\ell}): r\geq 1\} .$$
\item If $i_\ell=2\,i_{m+\ell}$ and $i_\ell\neq 0$, then $\mathcal L_+(D_H(x^i))$ contains
$$\{D_H(x^{i+r\cdot 1_\ell}): r\geq 2\}.$$
\end{enumerate}
Similar results hold when we swap $\ell$ for $m+\ell$ in $(1) $ and $(2)$. 
\end{lemma}

\begin{proof}
Let $r\geq 1$. Setting $j=(r+1)\cdot 1_\ell +1_{m+\ell}$, from \eqref{basishamilton}, we get
$$[D_H(x^i),D_H(x^j)]=\left(i_{m+\ell} (r+1) - i_{\ell}\right) D_H(x^{i+r\cdot 1_\ell}). $$  
So when $i_{m+\ell} (r+1) \neq i_{\ell}$, we see that
\begin{equation}\label{ham} 
\ell_+([D_H(x^i),D_H(x^j)])=D_H(x^{i+r\cdot 1_\ell})
\end{equation}
and furthermore \eqref{basishamilton} also yields that when $D_H(x^v)<D_H(x^i)$
\begin{equation}\label{ham2}
\ell_+([D_H(x^v),D_H(x^j)])<\ell_+([D_H(x^i),D_H(x^j)]).
\end{equation}
Hence, $D_H(x^{i+r\cdot 1_\ell})\in \mathcal L_+(D_H(x^i))$. 

\bigskip 

We now prove (1). Assume $i_\ell\neq 2i_{m+\ell}$. From what we have shown in \eqref{ham}, we may assume $i_{m+\ell} (r+1)= i_{\ell}$, letting $u=2\cdot 1_{\ell}+1_{m+\ell}$, from \eqref{basishamilton} we see that
$$[D_H(x^i), D_H(x^u)]=\left( 2\, i_{m+\ell} -i_\ell \right)D_H(x^{i+1_\ell}).$$
Note that the coefficient above is nonzero (as we are assuming $i_\ell\neq 2i_{m+\ell}$). It follows, using again \eqref{basishamilton} as in \eqref{ham2}, that $D_H(x^{i+1_\ell})\in \mathcal L_+(D_H(x^i))$. If $r=1$ we are done. On the other hand, if $r>1$, letting $j'=r\cdot 1_\ell +1_{m+\ell}$ we see that
$$[D_H(x^{i+1_\ell}), D_H(x^{j'})]=(i_{m+\ell}r-(i_{\ell}+1))D_H(x^{i+r\cdot 1_\ell}). $$
Since we are assuming $i_{m+\ell} (r+1)= i_{\ell}$, we get $i_{m+\ell}r\neq i_{\ell}+1$, so
$$\ell_+([D_H(x^{i+1_\ell}), D_H(x^{j'})])=D_H(x^{i+r\cdot 1_\ell})$$
 and again using \eqref{basishamilton} as in \eqref{ham2} we get $D_H(x^{i+r\cdot 1_\ell})\in \mathcal L_+(D_H(x^i))$, as desired. 
 
\bigskip
 
We now prove (2). Assume $i_\ell=2i_{m+\ell}$ and $i_\ell\neq 0$. These assumptions yield that $i_{m+\ell}(r+1)\neq i_\ell$. Hence from \eqref{ham} and \eqref{ham2} we get that $D_H(x^{i+r\cdot 1_\ell})\in \mathcal L_+(D_H(x^i))$, as desired.
\end{proof}

\begin{proposition}
With respect to the above ordering on $\MM$, the Hamiltonian Cartan algebra $\HH_n$ has no infinite leading-Dicksonian sequence. 
\end{proposition}
\begin{proof}
The proof is almost identical to the proofs of Proposition~\ref{forW} and \ref{forS}. We associate to each element $D_H(x^i)$ in the basis the nonzero tuple $i$ in $\mathbb N^n$. Now use Lemma~\ref{leadingH} to contradict Dickson's lemma. 
\end{proof}

\medskip

\subsection{Contact Cartan algebras}
Let $n=2m+1$ with $m$ a positive integer. The \emph{contact Cartan algebra} $\KK_n$ is the Lie subalgebra of $\WW_n$ given by elements of the form
\begin{align*}
D_K(p):=& \quad \sum_{\ell=1}^{m}\left(\partial_{m+\ell}(p)\partial_\ell \; -\; \partial_\ell(p)\partial_{m+\ell} \right) \; +\; \sum_{\ell=1}^{2m}x_\ell\partial_{n}(p)\partial_\ell \\
 & \quad + \left(2p \; -\; \sum_{\ell=1}^{2m} x_\ell\partial_\ell(p)\right)\partial_{n}
\end{align*}
for $p\in \kk[x_1,\dots,x_n]$. In fact $D_K:\kk[x_1,\dots,x_n]\to W_n$ as defined above is an injective 
 $\kk$-linear mapping. For $p,q\in \kk[x_1,\dots,x_n]$ we have $[D_K(p),D_K(q)]=D_K(\langle p,q\rangle)$ where 
\begin{equation}\label{forcontact}
\langle p,q\rangle=D_K(p)(q) \; -\; 2\partial_n(p)q .
\end{equation}
See \cite[\S1.3]{Hu}, for instance.

We prove that $\KK_n$ has a basis $\MM$ (of homogeneous elements) with an order (compatible with the natural grading) such that there is no infinite leading-Dicksonian sequence. Hence, by Theorem~\ref{basisPoisson}, the symmetric algebra $S(\KK_n)$ has ACC on radical Poisson ideals. 

\smallskip

Let $\MM=\{D_K(x^i): i\in \mathbb N^n\}\subset \KK_n$. The set $\MM$ is clearly a $\kk$-basis for $\KK_n$. From~\eqref{forcontact}, we see that
\begin{equation}\label{basiscontact}
[D_K(x^i),D_K(x^j)]=  D_K(D_K(x^i)(x^j) - 2 \partial_n(x^i) x^j).
\end{equation}
The inside term can be computed as
\begin{align*}
D_K(x^i)(x^j) - 2 \partial_n(x^i) x^j= &\;  \sum_{\ell=1}^{m}(i_{m+\ell}j_{\ell} \; -\; i_{\ell} j_{m+\ell})x^{i+j-1_\ell-1_{m+\ell}} \\
& \; + \sum_{\ell=1}^{2m} (i_nj_{\ell} \; -\; i_\ell j_n) x^{i+j-1_n}\\
& \; + 2(j_n\; -\;  i_n)x^{i+j-1_n}.
\end{align*}

For $i\in \mathbb N^n$, we set $|i|_K=\sum_{\ell=1}^{2m} i_\ell \; +\; 2i_n$. Recall that in the natural grading of $\KK_n$ the element $D_K(x^i)$ is homogeneous of degree $|i|_K-2$. Thus we equip $\MM$ with the following (total) order. Given $i,j\in \mathbb N^n$, we set 
$$D_K(x^i)< D_K(x^j) \quad \Leftrightarrow\quad (|i|_K,i_n,\dots,i_1)<_{\text{lex}}(|j|_K,j_n,\dots,j_1) .$$

We now prove a lemma which can be thought of as the $\HH_{n}$ analogue of  Lemmas~\ref{leadingW}, \ref{leadingS} and \ref{leadingH}.

\begin{lemma}\label{leadingK}
Let $i\in \mathbb N^n$. 
\begin{enumerate}
\item Let $1\leq \ell \leq m$. If $i_\ell \neq 2 \,i_{m+\ell}$, then $\mathcal L_+(D_K(x^i))$ contains
$$\{D_K(x^{i+r\cdot 1_\ell}): r\geq 1\}.$$
If $i_\ell=2\,i_{m+\ell}$ and $i_\ell\neq 0$, then $\mathcal L_+(D_K(x^i))$ contains
$$\{D_K(x^{i+r\cdot 1_\ell}): r\geq 2\}.$$
Similar results hold when we swap $\ell$ for $m+\ell$.  
\medskip
\item If $|i|\neq 2$, then $\mathcal L_+(D_K(x^i))$ contains
$$\{D_K(x^{i+r\cdot 1_n}): r\geq 1\}. $$
If $|i|=2$, then $\mathcal L_+(D_K(x^i))$ contains
$$\{D_K(x^{i+r\cdot 1_n}): r\geq 2\}. $$

\end{enumerate}
\end{lemma}

\begin{proof}
(1) Let $r\geq 1$. Setting $j=(r+1)\cdot 1_\ell +1_{m+\ell}$, from \eqref{basiscontact}, we get
\begin{multline*}
[D_K(x^i),D_K(x^j)]= \\
 \left(i_{m+\ell} (r+1) - i_{\ell}\right) D_K(x^{i+r\cdot 1_\ell}) \;  +\; i_n(r-1)\, D_K(x^{i+(r+1)\cdot 1_\ell + 1_{m+\ell}-1_n}). 
 \end{multline*}
Thus, due to the nature of the order in the basis $\MM$, when $i_{m+\ell} (r+1) \neq i_{\ell}$, we see that
$$
\ell_+([D_K(x^i),D_K(x^j)])=D_K(x^{i+r\cdot 1_\ell})
$$
and, using the two equalities above, one easily checks that for $D_K(x^v)<D_K(x^i)$ 
\begin{equation}\label{cont11}
\ell_+([D_K(x^v),D_K(x^j)])<\ell_+([D_K(x^i),D_K(x^j)]).
\end{equation}
Hence, $D_K(x^{i+r\cdot 1_\ell})\in \mathcal L_+(D_K(x^i))$.
 The rest of the argument follows the same lines as the proof of Lemma~\ref{leadingH}.

\bigskip 

(2) Let $r\geq 1$. Setting $j=(r+1)\cdot 1_n$, from \eqref{basiscontact}, we get
\begin{equation}\label{contact2}
[D_K(x^i),D_K(x^j)]= \; \left((r+1)(2-\sum_{\ell=1}^{2m}i_\ell)-2i_n\right)D_K(x^{i+r\cdot1_n}).
\end{equation}
Thus, when the coefficient on the right-hand-term is nonzero, we see that
\begin{equation}\label{contact} 
\ell_+([D_K(x^i),D_K(x^j)])=D_K(x^{i+r\cdot 1_n})
\end{equation}
and from these formulas one also readily checks that inequalities of the form \eqref{cont11} still hold. Hence, $D_K(x^{i+r\cdot 1_n})\in \mathcal L_+(D_K(x^i))$.

\medskip

Now let us assume $|i|\neq 2$, where recall that $|i|=\sum_{\ell=1}^{n}i_\ell$. From what we have shown in \eqref{contact}, we may assume that 
\begin{equation}\label{contact3}
(r+1)(2-\sum_{\ell=1}^{2m}i_\ell)=2i_n .
\end{equation}
Letting $u=2\cdot 1_{n}$, from \eqref{contact2} we see that
$$[D_H(x^i), D_H(x^u)]=\left(2(2-\sum_{\ell=1}^{n}i_\ell)\right)D_H(x^{i+1_n}) .$$
Note that the coefficient above is nonzero (as we are assuming $|i|\neq 2$). It then follows, after using these formulas to check that inequalities of the form \eqref{cont11} still holds, that $D_K(x^{i+1_n})\in \mathcal L_+(D_K(x^i))$. If $r=1$ we are done. On the other hand, if $r>1$, letting $j'=r\cdot 1_n$ we see that
$$[D_K(x^{i+1_n}), D_K(x^{j'})]=\left(r(2-\sum_{\ell=1}^{2m}i_\ell)-2(i_n+1)\right)D_K(x^{i+r\cdot 1_n}). $$
Since we are assuming \eqref{contact3}, we get that the above coefficient is nonzero, so
$$\ell_+([D_K(x^{i+1_n}), D_K(x^{j'})])=D_K(x^{i+r\cdot 1_n})$$
and again one checks using these formulas that inequalities of the form \eqref{cont11} hold. Hence, $D_K(x^{i+r\cdot 1_n})\in \mathcal L_+(D_K(x^i))$ as desired. 
 
\bigskip
 
Now let us assume $r\geq 2$ and  $|i|= 2$. This assumption yields that 
$$(r+1)(2-\sum_{\ell=1}^{2m}i_\ell)\neq 2i_n .$$ 
Hence from \eqref{contact} we get that $D_K(x^{i+r\cdot 1_n})\in \mathcal L_+(D_K(x^i))$, as desired.
\end{proof}

\medskip

\begin{proposition}
With respect to the above ordering on $\MM$, the contact Cartan algebra $\KK_n$ 
is Dicksonian.
\end{proposition}
\begin{proof}
The proof is almost identical to the proofs of Proposition~\ref{forW} and \ref{forS}. We associate to each element $D_K(x^i)$ in the basis the tuple $i$ in $\mathbb N^n$. Now use Lemma~\ref{leadingK} to contradict Dickson's lemma. 
\end{proof}

\medskip

\subsection{Loop algebras}
In this section we prove that loop algebras, their current subalgebras, and twisted loop algebras are Dicksonian.
Hence, by Theorem~\ref{basisPoisson}, their symmetric algebras have ACC on radical Poisson ideals. 

In this subsection, let $\kk$ be an algebraically closed field of characteristic zero.
\begin{theorem}\label{thm:loop}
Let  $\mf g$ be a simple finite-dimensional Lie algebra.  The loop algebra $\mf g[t, t^{-1}]$ and the current subalgebra $\mf g [t]$ have no infinite leading-Dicksonian sequences.
\end{theorem}

\begin{proof}
We give the proof for the loop algebra $\hat{\mf g} := \mf g[t,t^{-1}]$.  Recall that the Lie bracket on $\hat {\mf g}$ is given by $[g t^i, h t^j] = [g,h]t^{i+j}$, where $g, h \in \mf g$.

We first fix notation. 
See \cite{Humphreys} for terminology.
 Let $\mf h$ be a Cartan subalgebra 
 of $\mf g$.
Let $\Phi \subset \mf h^*$ be the set of roots of $\mf h$ acting on $\mf g$, and fix a set $\Delta \subset \Phi$ of simple roots; note that $\Delta$ is a basis for $\mf h^*$. Let $\Phi^+$ and $\Phi^-$ be respectively the set of positive and negative roots with respect to $\Delta$, so $\Phi = \Phi^- \sqcup \Phi^+$.  

Let $\mf g = \mf h \oplus \bigoplus_{\alpha \in \Phi} \mf g_\alpha$ be the root space decomposition of $\mf g$.  This gives a  grading of $\mf g$ by $\ZZ  \Delta$ since $[\mf g_\alpha, \mf g_\beta] \subseteq \mf g_{\alpha + \beta}$.  
By definition, if $x \in \mf g_\alpha$ and $h \in \mf h$, then $[h, x] = \alpha(h) x$.

We first define a (strict, total) order $\prec$ on $\Phi$.    Fix  an enumeration $ \{ \delta_1, \dots, \delta_r\}$ of $\Delta$.  Let $\alpha, \beta \in \Phi$, and write $\alpha = \sum_{i=1}^r a_i \delta_i$ and $\beta = \sum_{i=1}^r b_i \delta_i$.  
Let $\height(\alpha) := \sum a_i$. 
We say that 
$\alpha \prec \beta$ if and only if 
\beq \label{prec} 
\height(\alpha)  < \height(\beta), \mbox{ or } \height(\alpha) = \height(\beta) \mbox { and } (a_r, \dots, a_1) <_{\rm lex} (b_r, \dots, b_1).
\eeq 
(The reason for this definition is so that $\delta_1 \prec \delta_2 \prec \dots \prec \delta_r$.)
We make the convention that $\prec$ extends to an order on $\Phi \cup \{ 0\}$ by defining all elements of $\Phi^-$ to be $ \prec 0 $ and all elements of $ \Phi^+$ to be $\succ 0$.

For each $\alpha \in \Phi$, fix $0 \neq x_\alpha \in \mf g_\alpha$.   Let $\{ h_1, \dots, h_r\}$ be the basis of $\mf h$ dual to $\Delta$.  Let $\BB = \{ x_\alpha : \alpha \in \Phi\} \cup \{h_1, \dots, h_r\}$; as each vector space $\mf g_\alpha$ is one-dimensional, $\BB$ is a basis for $\mf g$.  
It is not always true that the bracket of two basis elements is a scalar multiple of a basis element because of our choice of basis for $\mf h$; however, if $\alpha+\beta \neq 0$ then $[x_\alpha, x_\beta]$ is a scalar multiple of a basis element.

Now define an order $<$ on $\BB$ by:
\begin{itemize}
\item $x_\alpha < x_\beta$ $\iff$ $\alpha \prec \beta $;
\item For all $\alpha \in \Phi^-$, $\beta \in  \Phi^+$, and $1 \leq i\leq r$, we have $x_{\alpha} < h_i < x_\beta$;
\item $h_1 < h_2 < \dots < h_r$.
\end{itemize}

We observe that, by simplicity of $\mf g$, there is some positive integer $K$ so that for all $x \in \BB$, the elements
\[ \{ [y_1, [y_2, \dots, [y_K, x] \dots ]] : y_1 , \dots, y_K \in \BB\}\]
span $\mf g$.

We now consider $\hat{\mf g}$, which is $\ZZ$-graded by degree in $t$, and also $\ZZ \times \ZZ \Delta$-graded, in the obvious way.  Define a triangular decomposition $\hat {\mf g} = \hat {\mf g}_- \oplus \hat {\mf g}_0 \oplus \hat {\mf g}_+$ via the $\ZZ$-grading, so $\hat {\mf g}_0 = \mf g$.
We extend $\BB$ to a basis  $\hat \BB := \{ x t^n : x \in \BB, n \in \ZZ\} $ of $\hat{\mf g}$, and define an order $<$ on $\hat \BB$ by $x t^m < y t^n$ $\iff$ $(m, x) <_{\rm lex} (n, y)$.
We refer to the corresponding order on  $\ZZ \times (\Phi\cup \{0\})$ as $ \prec $ as well.

Let $\hat \BB_+ = \hat \BB \cap \hat {\mf g}_+$ and $\hat \BB_- = \hat \BB \cap \hat {\mf g}_-$.  Let $\mc I_+ $ be the set of finite tuples from $\hat \BB_+$ and likewise let $\mc I_-$ be the set of finite tuples from $\hat \BB_-$.

Fix $M = x t^m \in \hat \BB$ (where $x\in \BB$ and $m \in \ZZ$), and let $\alpha$ be the $\ZZ \times \ZZ \Delta$-weight of $M$.
We analyze the sets $\mc L_-(M), \mc L_+(M)$.
Let $N = y t^n \in \hat{\BB}$ satisfy $N < M$.  Let $\beta$ be the $\ZZ \times \ZZ\Delta$-weight of $N$.  

Now if $x \not \in \mf h$ then $\hat {\mf g}_{\alpha}$ is one-dimensional, so  $N < M $ $\iff$  $\beta \prec \alpha$.  
Let ${\bf i} = (M_1, \dots, M_n ) \in \mc I_+$.  For $1 \leq i \leq n$, write $M_i \in \hat{\mf g}_{\alpha_i}$, where $\alpha_i \in \ZZ \times (\Phi\cup \{0\})$.
The $\ZZ \times \ZZ\Delta$-grading of $\hat {\mf g}$ means that $D_{\bf i}^+(M) \in \hat{\mf g}_{\alpha+ \sum \alpha_i}$.
Thus  $\beta + \sum \alpha_i \prec \alpha + \sum \alpha_i$ and if $D_{\bf i}^+(M), D_{\bf i}^+(N) \neq 0$ we have $ D_{\bf i}^+(N) < D_{\bf i}^+(M) $.
Thus $\mc L_+(M) = \{ D_{\bf i}^+(M) : {\bf i} \in \mc I_+, D_{\bf i}^+(M) \neq 0\}$.
By letting ${\bf i}$ be of the form $(y_1 t^{s}, y_2 t, \dots, y_K t)$ we see that $\mc L_+(M) \supseteq \mf \BB t^{\geq m+K}$.    Likewise $\mc L_-(M) \supseteq \mf \BB t^{\leq m-K}$.  

Now suppose that $x = h_j \in \mf h$, and consider ${\bf i}$ of the form $(y_1 t^s, y_2 t, \dots, y_K t, x_{\delta_j}t)$.  We may suppose that $D_{\bf i}^+(M) \neq 0$.
If $\beta = \alpha$ then $N = h_i t^m$ with $i < j$, and so $D_{\bf i}^+(N) = [x_{\delta_j}, h_i] = 0$.
If $\beta \prec \alpha$ then as before $D_{\bf i}^+(N) < D_{\bf i}^+(M)$.  In any case $D_{\bf i}^+(M) \in \mc L_+(M)$; letting the $y_i$ vary we see that $\mc L^+(M) \supseteq \mf \BB t^{> m+K}$.  Likewise, $\mc L_-(M) \supseteq \mf \BB t^{< m-K}$. 

By Remark~\ref{rem:newcond}, we see that $\hat{\mf g}$ has no infinite leading-Dicksonian sequences.
\end{proof}

A similar result holds for twisted loop algebras, which we now define.
Let $\mf g$ be a finite-dimensional simple Lie algebra, let $\sigma$ be an automorphism of $\mf g$ of order $m$, and let $\eta$ be a primitive $m$-th root of unity. 
Each eigenvalue of $\sigma$ has the form $\eta^j$ for $ j \in \ZZ_m$, and this gives a $\ZZ_m$-grading of $\mf g$, which we write $\mf g = \bigoplus_{j \in \ZZ_m} \mf g_j$.
If $ j\in \ZZ$ let $\overline{j} = j \mod m$.
The twisted loop algebra $L(\mf g, \sigma, m)$ 
 is the $\ZZ$-graded Lie subalgebra of $\hat{\mf g}$ with  $L_j = \mf g_{\overline j} t^j$.

\begin{theorem} Let $\mf g$ be a finite-dimensional simple Lie algebra, let $\sigma$ be an automorphism of $\mf g$ of order $m$, and let $\eta$ be a primitive $m$-th root of unity.  Then $L := L(\mf g, \sigma, m)$ is Dicksonian.
\end{theorem}
\begin{proof}
The proof is similar to the proof of Theorem~\ref{thm:loop}, and we omit the details.
\end{proof}

\medskip

\subsection{Simple graded Lie algebras of polynomial growth}

We now prove Corollary~\ref{simplecase}. We restate it for the reader's convenience.

\medskip

\begin{corollary}\label{cor:simplerestate}
Let $\kk$ be a field of characteristic zero and $\gg$ a Lie $\kk$-algebra. If $\kk^{\text{alg}}\otimes_\kk\gg$ is a simple graded Lie $\kk^{\text{alg}}$-algebra of polynomial growth, then the symmetric algebra $S(\gg)$ has ACC on radical Poisson ideals.
\end{corollary}

When $\kk$ is algebraically closed, simple graded Lie algebras of polynomial growth have been classified by Mathieu as follows:

\begin{theorem}[Mathieu's Classification \cite{mathieu}]\label{mathieuclass}
Assume $\kk$ is an algebraically closed field of characteristic zero. Let $\gg$ be a $\mathbb Z$-graded simple Lie $\kk$-algebra of polynomial growth. Then, $\gg$ is one of the following
\begin{enumerate}
\item a finite-dimensional simple Lie algebra; or
\item a (twisted or untwisted) loop algebra; or
\item a Cartan type algebra of the form $\WW_n, \SS_n, \HH_n$, or $ \KK_n$; or
\item the Witt algebra.
\end{enumerate}
\end{theorem}

Due to the above classification, 
{
the proof of Corollary~\ref{cor:simplerestate} will follow as a} consequence of the previous results in this section together with the following general result 
. Recall that given a Poisson $\kk$-algebra $(A,\{-,-\})$ and a field extension $L/k$ we can equip $L\otimes_\kk A$ with the canonical Poisson bracket where the Poisson structure on $L$ is trivial and so $A$ becomes naturally a Poisson $L$-algebra. Namely,
$$\{a\otimes b,c\otimes d\}=ac\otimes \{b,d\},\quad \text{ for } a,c\in L \text{ and }b,d\in A.$$
Below we assume $L\otimes_\kk A$ is equipped with this Poisson bracket.

\begin{lemma}\label{extscalars}
Let $(A,\{-,-\})$ be a Poisson algebra over a field $\kk$ of characteristic zero and let $L$ be a field extension of $\kk$. If the Poisson $L$-algebra $L\otimes_\kk A$ has ACC on radical Poisson ideals then the same holds in $A$.
\end{lemma}
\begin{proof} 
First note that, due to the nature of the Poisson bracket of $L\otimes_\kk A$, if $I$ is a Poisson ideal of $A$ then $L\otimes_\kk I$ is a Poisson ideal of $L\otimes_\kk A$. It suffices to show that if $I_1$ and $I_2$ are radical Poisson ideals of $A$ with $I_1$ properly contained in $I_2$, then the same holds for $L\otimes_\kk I_1$ and $L\otimes_\kk I_2$. Since $\kk$ is of characteristic zero, the ideals $L\otimes_\kk I_i$ are radical (recall that being a reduced ring is preserved under base change over perfect fields). 
By faithful flatness ($\kk$ being a field), $L\otimes_\kk I_1$ is properly contained in $L\otimes_\kk I_2$. The result follows.
\end{proof}

{
\begin{proof}[Proof of Corollary \ref{cor:simplerestate}]
Let $\gg'=\kk^{\text{alg}}\otimes_{\kk}\gg$. Namely, $\gg'$ is the canonical base change of $\gg$ to a Lie $\kk^{\text{alg}}$-algebra. Thus, the Lie bracket on $\gg'$ is the unique $\kk^{\text{alg}}$-bilinear map extending that on $\gg$.  Our assumption states that $\gg'$ is a simple graded Lie $\kk^{\text{alg}}$-algebra of polynomial growth.


By Theorem \ref{mathieuclass}, $\gg'$ is one of the Lie algebras from the list (1)-(4) in the statement of that theorem. In case (1), $\gg'$ is finite-dimensional and so $S(\gg')$ is outright noetherian (by Hilbert's basis theorem). For cases (2), (3) and (4), we have already established that $\gg'$ is Dicksonian (in the work done previously in this section), and so, by Theorem~\ref{basisPoisson}, $S(\gg')$ has ACC on radical Poisson ideals.

Finally, since 
$$S(\gg')=S(\kk^{\text{alg}}\otimes_{\kk} \gg)\cong \kk^{\text{alg}}\otimes_{\kk}S(\gg),$$
by Lemma \ref{extscalars}, $S(\gg)$ has ACC on radical Poisson ideals, as desired. 
\end{proof}
}

\section{On the Poisson Dixmier-Moeglin Equivalence}\label{PoiDME}

In this final section we make some remarks on the Poisson Dixmier-Moeglin equivalence in the context of Poisson algebras of countable dimension (for instance, those that are finitely Poisson-generated). Let $\kk$ be a field (in arbitrary characteristic unless stated otherwise) and $(A,\{-,-\})$ a Poisson $\kk$-algebra. Let $P$ be a prime Poisson ideal of $A$. We recall that $P$ is Poisson-primitive if there is a maximal ideal $M$ of $A$ such that $P$ is the largest Poisson ideal contained in $M$ (i.e., $P$ is the Poisson core of $M$). On the other hand, $P$ is said to be Poisson-locally closed if $P$ is a locally closed point in the Poisson spectrum of $A$; and $P$ is said to be Poisson-rational if the Poisson centre of the field of fractions of $A/P$ is algebraic over $\kk$ 
{
(recall that the Poisson centre consists of those elements $a$ such that $\{a,-\}$ is the trivial derivation)}.

An algebra is said to satisfy the Poisson Dixmier-Moeglin equivalence (or PDME) if the notions of Poisson-primitive, Poisson-locally closed, and Poisson-rational, coincide. We refer the reader to the introduction of \cite{BLLSM} for further details and recent developments on the subject. Under mild assumptions, we prove in Theorem~\ref{implications} that one has the following implications:
$$\text{Poisson-locally closed} \quad \implies \quad \text{Poisson-primitive} \quad \iff \quad \text{Poisson-rational} . $$
To prove that Poisson-rational implies Poisson-primitive we will use the following result from \cite[Lemma 3.1]{BLLSM}.

\begin{lemma}\label{newelement}
 Let $\kk$ be a field and $A$ an integral and commutative $\kk$-algebra equipped with a collection of $\kk$-linear derivations $(\delta_j)_{i\in J}$.
Suppose that there is a finite-dimensional $\kk$-vector subspace $V$ of $A$ and a set $\mathcal{S}$ of ideals satisfying:
\begin{itemize}
\item[(i)]
$\delta_j(I)\subseteq I$ for all $j\in J$ and $I\in\mathcal S$,
\item[(ii)]
$\bigcap\mathcal{S} = (0)$, and
\item[(iii)]
 $V\cap I\neq (0)$ for all $I\in \mathcal{S}$.
\end{itemize}
Then there exists $f$ in Frac$(A)\setminus \kk$ with $\delta_j(f)=0$ for all $j\in J$.
\end{lemma}

 \begin{remark}\label{remarknew} \
 \begin{enumerate}
 \item We point out that in \cite[Lemma 3.1]{BLLSM} the above statement appears in the case when $J$ is finite (namely, a finite collection of derivations). However, the proof does not use finiteness of $J$ and could have been stated there for general indexing set $J$.
 \item We also note that the conclusion of the lemma can be strengthen to find $f$ not algebraic over $\kk$. Indeed, if we let $\kk^{\alg}$ be the relative algebraic closure of $\kk$ in $A$, then we can view $A$ as a $\kk^{\alg}$-algebra and the derivations $(\delta_j)_{j\in J}$ are $\kk^{\alg}$-linear. Letting $V'$ be the span of $V$ over $\kk^{\alg}$, we see that the conditions in the lemma still hold when replacing $V$ for $V'$. Thus, we find $f\in Frac(A)\setminus \kk^{\alg}$ as desired.
 \end{enumerate}
 \end{remark}

\begin{theorem}\label{implications}
Let $\kk$ be an uncountable field and $(A,\{-,-\})$ a Poisson $\kk$-algebra that has countable dimension over $\kk$ (for example, when $A$ is finitely Poisson-generated over $\kk$). Then, for a prime Poisson ideal
$$\text{Poisson-locally closed} \quad \implies \quad \text{Poisson-primitive} \quad \iff \quad \text{Poisson-rational} . $$
\end{theorem}
\begin{proof}
By \cite{amitsur}, the assumptions on $\kk$ and $A$ yield that $A$ is a Jacobson ring. Hence, \cite[Proposition 1.7(i)]{oh} yields that Poisson-locally closed implies Poisson-primitive. 

\medskip

On the other hand, for $P$ a Poisson-primitive ideal of $A$ with corresponding maximal ideal $M$, in the proof of \cite[Proposition 1.10]{oh} an injective morphism from the Poisson centre of the field of fractions of $A/P$ to $End_A(A/M)$ is constructed. By our assumptions on $\kk$ and $A$, $A$ satisfies the Nullstellensatz; in particular, $End_A(A/M)$ is algebraic over $\kk$. Thus, $P$ is Poisson-rational.

\medskip

Finally, to show that Poisson-rational implies Poisson-primitive we adapt the argument from \cite[Theorem 3.2]{BLLSM} \footnote{We thank Alexey Petukhov for pointing out the necessary adaptations.}. Without loss of generality, we may assume that $P=(0)$ is Poisson-rational (as we may replace $A$ for $A/P$ if necessary). Let $\mathcal S$ be the collection of all nonzero proper Poisson ideals of $A$. 
{
We may assume that $\mathcal S$ is non empty, as otherwise $(0)$ is the Poisson core of any maximal ideal of $A$ and hence Poisson-primitive}. Now, since $A$ has countable dimension, there is a chain of finite-dimensional $\kk$-vector subspaces $V_1\subset V_2\subset \cdots$ such that 
$$A=\bigcup_{n\geq 1} V_n .$$
Set
$$\mathcal S_n=\{Q\in \mathcal S: Q\cap V_n\neq (0) \} .$$
Note that $\mathcal S=\bigcup_n\mathcal S_n$. We claim that each 
{
intersection $\bigcap\mathcal S_n$ is nontrivial}
{
(recall that when $\mathcal S_n$ is empty the intersection $\bigcap\mathcal S_n$ is generally defined to be $\mathcal S$)}. Towards a contradiction, assume $\bigcap \mathcal S_n=(0)$ 
{
for some $n$}. Let $(a_j)_{j\in J}$ be generators of $A$ as a $\kk$-algebra and let $\delta_j=\{a_j,-\}$ be the Hamiltonian derivation associated to $a_j$ for each $j\in J$. As the ideals in $\mathcal S_n$ are Poisson, they are also differential with respect to $(\delta_j)_{j\in J}$. We can now apply Lemma~\ref{newelement}, also see Remark~\ref{remarknew}(2), to get $f\in Frac(A)\setminus \kk^{\alg}$ such that $\delta_{j}(f)=0$ for all $j\in J$. But the latter equalities imply that $f$ is in the Poisson centre of the fraction field of $A$, contradicting Poisson-rationality of $(0)$. 

Now let $L_n=\bigcap \mathcal S_n$ for $n\geq 1$. We have shown that $L_n\neq (0)$, so let $f_n$ be a nonzero element in $L_n$. If we let $T$ be the (countable) multiplicatively closed set generated by the $f_n$'s, we see that the localization $B:=T^{-1}A$ is a countably generated $\kk$-algebra. Hence $B$ satisfies the Nullstellensatz (as $\kk$ is uncountable). If we let $I$ be any maximal ideal of $B$ and $J:=I\cap A$, then $A/J$ embeds into $B/I$ and the latter is an algebraic extension of $\kk$ (as $B$ satisfies the Nullstellensatz). Hence, $A/J$ is also an algebraic extension of $\kk$, thus a  field, and so $J$ is a maximal ideal of $A$. By construction, $J$ does not contain any element in $\mathcal S$ (since $\mathcal S=\bigcup_n\mathcal S_n$), and so $(0)$ is the largest Poisson ideal in $J$. In other words, $(0)$ is Poisson-primitive as desired. 
\end{proof}

Thus, in the cases of interest (finitely Poisson-generated complex Poisson algebras, for instance), the PDME reduces to showing that Poisson-primitive implies Poisson-locally closed. It is shown in \cite{BLLSM} that there are finitely generated Poisson algebras that do not satisfy the PDME (namely, have Poisson-primitive ideals that are not Poisson-locally closed). However, these examples are far from being symmetric algebras, and so one can ask whether the PDME holds for symmetric algebras $S(\gg)$ that are finitely Poisson-generated (and have ACC on radical Poisson ideals). When the Lie algebra $\gg$ is finite-dimensional and $\kk$ has characteristic zero, $S(\gg)$ does satisfy the PDME; this appears in \cite[Theorem 5.7]{LLS} (and can be thought of as the Poisson analogue of the seminal work of Dixmier and Moeglin~\cite{Dixmier,Moeglin} showing that the enveloping algebra $U(\gg)$ satisfies the classical DME). But in general the answer is no; for instance, when $\gg$ is the positive Witt algebra $W_+$ we have:

\begin{lemma}\label{easy}
Assume $\kk$ is of characteristic zero and let $W_+$ be the positive Witt algebra. In $S(W_+)$ the zero ideal $(0)$ is Poisson-rational but not Poisson-locally closed. 
\end{lemma}
\begin{proof}
Recall that 
$$W_+=\text{span}_{\kk}(e_i:i\geq 1).$$
If we let $P_i$ be the Poisson ideal of $S(W_+)$ generated by $e_i$, we see that each $P_i$ is a nonzero prime Poisson-ideal and $\bigcap_i P_i=(0)$. Thus $(0)$ is not Poisson-locally closed.

On the other hand, let $f/g$ be in a nonzero element in the centre of $Frac(S(W_+))$. We may assume that $f$ and $g$ have no common factors (recall that $S(W_+)$ is a UFD). Let $e_n$ be larger than any $e_i$ appearing in $f$ and $g$. As $f/g$ is in the centre, we have $\{f/g,e_n\}=0$ and this yields 
$$\{f,e_n\}g=f\{g,e_n\}.$$
As $f$ and $g$ have no common factors, by the choice of $e_n$, we must have $\{f,e_n\}=0$ and $\{g,e_n\}=0$. But this can only happen if $f$ and $g$ are in $\kk$. Thus $(0)$ is Poisson-rational. We note that this is also proved in \cite[Lemma~3.3]{Bois}.
\end{proof}

\begin{remark}
By Theorem~\ref{implications}, if $\kk$ is uncountable of characteristic zero, then $(0)$ is a Poisson primitive ideal of $S(W_+)$.
\end{remark}

Nonetheless, 
in the presence of ACC on radical Poisson ideals, one gets close to a PDME, as we point out in the lemma below. 

\begin{remark}\label{rem:pseudoPDME}
If $A$ is a Poisson algebra with ACC on radical Poisson ideals, then a prime Poisson-ideal $P$ is Poisson-locally closed iff there are finitely many (but at least one) prime Poisson ideals of Poisson-height one in $A/P$. Indeed, if $P$ is Poisson-locally closed then the intersection of all nonzero prime Poisson ideals in $A/P$ is nonzero; as this intersection is a nonzero radical Poisson ideal, the ACC and Theorem~\ref{conservativePoisson} yield that it must be a finite intersection of nonzero prime Poisson ideals of Poisson-height one (i.e., the Poisson components). Thus there are finitely many (and at least one) prime Poisson ideals of Poisson-height one in $A/P$. The other direction is clear. 
\end{remark}

\begin{lemma}\label{PDME}
Assume $A$ is a Poisson $\kk$-algebra of countable dimension over $\kk$ that has ACC on radical Poisson ideals. Let $P$ be a prime Poisson ideal of $A$. If $P$ is Poisson-rational then there are at most countably many prime Poisson ideals of Poisson-height one in $A/P$.
\end{lemma}
\begin{proof}
We may assume that $P=(0)$. Let $S$ denote the set of nonzero prime Poisson ideals of $A$ of Poisson-height one. We must show that $S$ is countable. Let $V_1\subset V_2\subset\cdots$ be a chain of finite-dimensional $\kk$-subspaces such that $A=\bigcup_i V_i$. Set
$$S_i=\{Q\in S: Q\cap V_i\neq (0)\}. $$
Note that $S=\bigcup_i S_i$, and so it suffices to show that $S_i$ is countable. We actually show each $S_i$ is finite. The same argument as in the proof of Theorem~\ref{implications} shows that the intersection $\bigcap S_i$ is nonzero. As the latter is a nonzero radical Poisson ideal, the ACC assumption and Theorem~\ref{conservativePoisson}(ii), yield that it must be a finite intersection of nonzero prime Poisson ideals of Poisson-height one (i.e., the Poisson components). This finite collection is in fact formed by all the elements in $S_i$, and so $S_i$ is finite.
\end{proof} 

As a consequence of Lemmas  \ref{easy} and \ref{PDME}, we see that in $S(W_+)$ there are either countably infinite many prime Poisson-ideals of Poisson-height one or there are none. We currently do not know the answer to this, so we leave it as an open problem:

\begin{question}
Are there any prime Poisson ideals in $S(W_+)$ of Poisson-height one (equivalently, of finite Poisson-height)?
\end{question}

The situation is quite different for nonzero Poisson ideals of $S(W_+)$. More precisely, we conclude by showing that if $P$ is a nonzero Poisson-rational ideal of $S(W_+)$ then there are only finitely many (and at least one) prime Poisson ideals of height one in $S(W_+)/P$ (the latter height being the classical algebraic one). 

\begin{proposition}\label{prop6.8}
Assume $\kk$ is of characteristic zero. Let $P$ be a nonzero prime Poisson ideal in $S(W_+)$. If $P$ is Poisson-rational, then there are finitely many (and at least one) prime Poisson ideals in $S(W_+)/P$ of (algebraic) height one.  
\end{proposition}
\begin{proof}
 Note that for any $e_i,e_j\in W_+$, the $\kk$-span of $\LL_+(e_i)\cup\LL_-(e_j)$ has finite codimension in $W_+$. Thus, by Corollary~\ref{cor:fingen}, there is $h\in S(W_+)\setminus P$ such that the localisation $(S(W_+)/P)_{h}$ is a finitely generated $\kk$-algebra.   Let $S_1$ be the collection of prime Poisson-ideals in $S(W_+)/P$ of height one that contain $h$ and $S_2$ those that do not contain $h$. Since $\bigcap S_1$ is nonzero (as it contains $h$), the ACC on radical Poisson ideals implies that this intersection has finitely many  Poisson components. Thus $S_1$ is finite. On the other hand, each element of $S_2$ yields a prime-Poisson ideal of height one in the finitely generated algebra $(S(W_+)/P)_{h}$. By \cite[Theorem~7.1]{BLLSM}, the fact that $P$ is Poisson-rational implies that this latter collection is finite, and so $S_2$ is finite. The result follows. 
\end{proof}

\begin{remark}\label{rem:pseudoPDME2} \
\begin{enumerate}
\item [(i)] We note that the proof of Proposition~\ref{prop6.8} works for any $\mf g$ with the property 
that the basis $\MM$ of $\mf g$ satisfies \eqref{cofinite}, 
as then we can invoke Corollary~\ref{cor:fingen}. In particular, Proposition~\ref{prop6.8} applies also to the symmetric algebra of the Witt algebra $W$, the (first) Cartan algebra $\WW_1$, and any (twisted or untwisted) loop algebra.
\item [(ii)] We also note that 
Proposition~\ref{prop6.8} does not imply that any nonzero Poisson-rational $P$ must be Poisson-locally closed (as it only refers to algebraic height, not Poisson-height). In fact, we currently do not know whether this is the case or not.\end{enumerate}
\end{remark}

\

\bibliographystyle{amsalpha}

\end{document}